\documentclass[a4paper, reqno]{amsart}
 \pdfoutput=1

\usepackage[utf8]{inputenc}
\usepackage{amsmath}
\usepackage{mathtools}
\usepackage{mathrsfs} 
\usepackage{amsfonts}
\usepackage{amssymb}
\usepackage{amsthm}
\usepackage{wasysym}
\usepackage{indentfirst}
\usepackage{graphicx}
\usepackage{enumitem}
\usepackage{xfrac}
\usepackage{esint}
\usepackage{comment}
\usepackage{placeins}

\usepackage{multirow}
\usepackage{calc}
\newcommand\Tstrut{\rule{0pt}{2.6ex}}         
\newcommand\Bstrut{\rule[-0.9ex]{0pt}{0pt}}   

\usepackage[backend=biber,maxbibnames=5,maxalphanames=5,style=alphabetic,bibencoding=utf8,giveninits,url=false,isbn=false,doi=true,eprint=true]{biblatex}
\renewbibmacro{in:}{}
\AtBeginBibliography{\small}

\usepackage{xcolor}
\definecolor{cobalt}{rgb}{0.0, 0.28, 0.67}
\definecolor{darkcerulean}{rgb}{0.03, 0.27, 0.49}

\usepackage{a4wide}

\usepackage[font=footnotesize,margin=0.2cm,skip=3pt]{caption}
\usepackage{subcaption}
\graphicspath{{./images/}}

\def\BS{{X}}

\def\del{\partial}
\def\eps{\varepsilon}
\def\zh{{z_{\tau}}}

\def\bilJ{{j}}
\def\bilR{{r}}
\def\bilu{{b}}
\def\proj{\Pi}

\def\softd{{\leavevmode\setbox1=\hbox{d}%
    \hbox to 1.05\wd1{d\kern-0.4ex{\char039}\hss}}}

\newcommand{\rP}{\mathcal{P}}

\newcommand{\cH}{\mathcal{H}}

\newcommand{\rV}{\mathbb{V}}

\providecommand{\normtmp}[2]{{#1\lVert #2 #1\rVert}}
 \providecommand{\norm}[1]{\normtmp{}{#1}}

  \providecommand{\skptmp}[3]{{\ensuremath{#1\langle {#2}, {#3} #1\rangle}}}
\providecommand{\skp}[2]{\skptmp{}{#1}{#2}}

\providecommand{\abstmp}[2]{{#1\lvert #2 #1\rvert}}
\providecommand{\abs}[1]{\abstmp{}{#1}}
\providecommand{\bigabs}[1]{\abstmp{\big}{#1}}

\providecommand{\setR}{\mathbb{R}}
\providecommand{\esssup}{\mathop{\mathrm{ess\,sup}}}

\def\div{\mathrm{div}}
\newcommand{\ds}{\,\mathrm{d}s}

\def\identity{\mathrm{Id}}
\def\tauref{\tau_{\textup{fine}}}
\def\zs{z_{\textup{s}}}

\usepackage{circuitikz}

\usepackage{todonotes}

\newcommand{\dt}{\,\mathrm{d}t}
\newcommand{\dx}{\,\mathrm{d}x}
\newcommand{\dsigma}{\,\mathrm{d}\sigma}

\newcommand{\tp}{^\mathsf{T}}

\newcommand{\quadnodes}{s_Q}
\newcommand{\projnodes}{s_{\Pi}}

\usepackage[normalem]{ulem}
\newcommand{\AKdeleted}{%
 \bgroup\markoverwith
  {\textcolor{red}{\pgfsetfillopacity{0.2}\rule[-0.5ex]{2pt}{10pt}\pgfsetfillopacity{1}}%
   \textcolor{red}{\llap{\rule[0.4ex]{2pt}{1pt}}}%
  }%
  \ULon}

\usepackage{pgf}
\usepackage{hyperref}
  \hypersetup{
  	breaklinks=true,
    citecolor=darkcerulean,
    linkcolor=darkcerulean,
    urlcolor=darkcerulean,
    pdfborder={0 0 0},
    pdfborderstyle={/S/U/W 0},
    colorlinks=true,
            }
            
\newtheorem{theorem}{Theorem}
\numberwithin{theorem}{section}
\theoremstyle{definition}

\newtheorem{example}[theorem]{Example}
\newtheorem{lemma}[theorem]{Lemma}

\newtheorem{proposition}[theorem]{Proposition}
\newtheorem{scheme}[theorem]{Scheme}
\newtheorem{assumption}[theorem]{Assumption}
\theoremstyle{remark}
\newtheorem{remark}[theorem]{Remark}

\title[Energy-consistent time discretization of port-Hamiltonian systems]{Energy-consistent Petrov--Galerkin time discretization of port-Hamiltonian systems}
\date{\today}

\makeatletter
\@namedef{subjclassname@2020}{%
	\textup{2020} Mathematics Subject Classification}
\makeatother

\author[J.~Giesselmann]{Jan Giesselmann}
\author[A.~Karsai]{Attila Karsai}
\author[T.~Tscherpel]{Tabea Tscherpel}

\address[J.~Giesselmann, T. Tscherpel]{Department of Mathematics, Technische Universität Darmstadt, Dolivostr.~15, 64293 Darmstadt,	Germany}
\address[A.~Karsai]{Institute of Mathematics, Technische Universität Berlin, Str.~des~17.~Juni 136, 10623 Berlin,	Germany}

\email{giesselmann@mathematik.tu-darmstadt.de}
\email{karsai@math.tu-berlin.de}
\email{tscherpel@mathematik.tu-darmstadt.de}

\begin{document}

\begin{abstract}
For a general class of nonlinear port-Hamiltonian systems we develop a high-order time discretization scheme with certain structure preservation properties. 
The finite or infinite-dimensional system under consideration possesses a Hamiltonian function, which represents an energy in the system and is conserved or dissipated along solutions.  
For infinite-dimensional systems this structure is preserved under suitable Galerkin discretization in space. 
The numerical scheme is energy-consistent in the sense that the Hamiltonian of the approximate solutions at time grid points  behaves accordingly.
This structure preservation property is achieved by specific design of a continuous Petrov--Galerkin (cPG) method in time.  
It coincides with standard cPG methods in special cases, in which the latter are energy-consistent. 
Examples of port-Hamiltonian ODEs and PDEs are presented to visualize the framework. 
In numerical experiments the energy consistency is verified and the  convergence behavior is investigated. 
\end{abstract}

\subjclass[2020]{
	35K55, 
	37L65, 
	37K58, 
	65J08, 
	65J15, 
	65L60, 
	65M60, 
	65P10
	}

\keywords{port-Hamiltonian system, Hamiltonian system, gradient system, energy conserving, structure preservation,  Petrov-Galerkin}

\maketitle

\section{Introduction}\label{sec:intro}
\noindent 
The framework of port-Hamiltonian systems allows to model complex physical systems of ordinary and partial differential equations that obey inherent energy conservation and dissipation principles. 
Applications arise for example in mechanics, in electronics, and in energy systems, see e.g.,~\cite[Ch.~6]{vanderschaft17-l2gain} and \cite{VanderSchaft2014,MehrmannUnger2023}. 
Also various complex fluid systems~\cite{Grmela1997} feature a Hamiltonian structure. 

\subsubsection*{Port-Hamiltonian systems}
Most infinite-dimensional port-Hamiltonian systems  are of the form
	\begin{align}\label{eq:ph-intro}
	C(z(t)) \partial_t z(t)
		= J(\mathcal{H}'(z(t))) 
		- R (\mathcal{H}'(z(t)))
		+ B(t,\mathcal{H}'(z(t))) \quad \text{ for } t  \geq t_0 ,
	\end{align}
with time-dependent state function $z$,	for possibly nonlinear operators $C, J, R, B$ (sometimes referred to as state dependent operators) with $C(z)$ acting on $\partial_t z$,  subject to initial conditions on $z$. 
In this work we focus on the subclass with $C(z) = \identity$.
Then, the system~\eqref{eq:ph-intro} reads
	\begin{align}\label{eq:ph-intro-2}
 \partial_t z(t)
	&= J(\mathcal{H}'(z(t))) 
	- R (\mathcal{H}'(z(t)))
	+ B(t,\mathcal{H}'(z(t))) \quad \text{ for } t  \geq t_0. 
\end{align}
The precise formulation is contained in Section~\ref{sec:prob-general}. 
In this formulation the operator $J$ describes conservative effects and the operator $R$ describes dissipative effects. 
The term $B$ contains all interactions with the environment, e.g., controls. 
There is a Hamiltonian $\mathcal{H}$ associated with this system, which
 is non-increasing along sufficiently smooth solutions $z$ if $B \equiv 0$, see Lemma~\ref{lem:pow-bal}. 
If additionally there is no dissipation in the system, i.e., $R \equiv 0$, then the Hamiltonian is conserved.
Depending on the properties of the operators, the formulations~\eqref{eq:ph-intro} and~\eqref{eq:ph-intro-2} contain port-Hamiltonian ODEs, PDEs and for degenerate $C$ the system~\eqref{eq:ph-intro} includes also differential algebraic equations (DAEs).

\subsubsection*{Structure preservation} 
To ensure robust computations, structure preservation properties of the numerical schemes are of particular interest. 
Especially for Hamiltonian ODEs structure-preserving schemes have a long history, see, e.g.,~\cite{Hairer2010} for an overview. 
In the context of (port)-Hamiltonian systems various different notions of \textit{structure preservation} are available, such as symplecticity~\cite{Hairer2010}, the conservation in time of the Hamiltonian by the approximate solutions~\cite{hairer10-energy,EggerHabrichShashkov2021,morandin24-modeling}, or the preservation of the (port)-Hamiltonian structure under spatial discretization~\cite{EggerHabrichShashkov2021,morandin24-modeling}. 
Here we focus on the latter two types, in the following referred to as \emph{energy consistency} and as \emph{structure-preserving space discretization}, respectively. 
We call a scheme energy-consistent if 
for any port-Hamiltonian system the energy of solutions to the scheme has the same behavior at all time grid points as the energy of the exact solution; see Proposition~\ref{prop:discr-energy} for the precise notion.
In particular, for systems without control ($B \equiv 0$) the Hamiltonian is (exactly) non-increasing, and if additionally there is no dissipation ($R \equiv 0$), then the Hamiltonian is exactly preserved, cf.~also~\cite[Def.~III.2]{celledoni17-energy} and \cite[Def.~1]{kotyczka19-discrete}. 
This is a stronger notion of energy consistency than the one used, e.g., in~\cite[Sec. IIIC]{mehrmann19-structure}, for a comparison see Remark~\ref{rmk:comp-energyconsist} below. 
Note that due to the celebrated Ge--Marsden Theorem~\cite{Ge1988} energy-consistent schemes (as introduced before) with fixed time step size, cannot be symplectic for Hamiltonian systems. 
For this reason our scheme is not symplectic. 
We call a space discretization of a port-Hamiltonian system \eqref{eq:ph-intro-2} \emph{structure-preserving}, if the space discrete system has again the form~\eqref{eq:ph-intro-2}. 
Numerical schemes with these properties are presented in~\cite{EggerHabrichShashkov2021, morandin24-modeling}. 

 \subsubsection*{Available methods}
Exact energy preservation for Hamiltonian and energy dissipation for gradient systems has attracted a lot of attention, and there is a wide range of methods available including some of high order. 
 
\emph{Discrete gradient methods} date back to~\cite{gonzalez96-time,McLachlanQuispelRobidoux1999} and usually are exactly energy-preserving for Hamiltonian systems and of second order. 
By now, also high-order generalizations are available, among  others~\cite{eidnes22-order}, see also the references therein. 
Recently, in~\cite{schulze23-structure2} a second order discrete gradient method has been developed for port-Hamiltonian systems with state-dependent operator $C(z)$.

Another class of methods are the averaged vector field collocation methods, also referred to as energy-preserving \emph{collocation methods}~\cite{hairer10-energy,cohen11-linear,Hairer2014}. 
They are exactly energy-preserving for Hamiltonian systems, and energy-dissipating for gradient systems. 
Furthermore, they have been applied to port-Hamiltonian systems with $C(z) = \identity$ and to gradient flows on Riemannian manifolds in~\cite{celledoni17-energy,Celledoni2018}. 
Note that they are of arbitrary order and they are energy-consistent. 

A structure preserving DG space discretization was presented in \cite{Jackaman2019} for the linearized KdV equation. 
In combination with the Crank--Nicolson scheme in time it is proven to be energy-consistent. 
In~\cite{JackamanPryer2021} this approach was extended to  energy-consistent schemes of lowest order in time for some dispersive Hamiltonian systems. 

Further methods include the \emph{continuous} and  \emph{discontinuous Petrov--Galerkin methods} (cPG and dPG for short); see, e.g., \cite[Sec.~69, 70]{Ern2021} for their formulation for general evolution equations. 
DPG methods are energy-dissipating rather than energy-preserving for Hamiltonian systems with convex Hamiltonian.
On the other hand, cPG methods have inherent energy consistency properties as shown in~\cite{French1990} for gradient flows (ODEs and PDEs).
It can easily be verified that standard cPG methods are energy-consistent for Hamiltonian systems with quadratic Hamiltonian and linear operator~$J$. 
For non-quadratic Hamiltonian or nonlinear operator~$J$, this is not true in general. 
Indeed, combining the lowest order cPG with the midpoint rule yields the Crank--Nicolson method. 
Hence, it is evident that energy consistency is not available in general for non-quadratic Hamiltonian.  
The cPG methods have been well-investigated, for example in~\cite{Aziz1989} for the heat equations. 
Therein, the authors use Gauß quadrature for the numerical approximation of the nonlinear terms and prove superconvergence at the time grid points. 
In~\cite{Schieweck2010} the cPG method (referred to as dPG method therein) was proven to dissipate energy for linear ODEs and gradient flows.  
It should be noted, that for linear systems with $B \equiv 0$ it reduces to a collocation method. 
In combination with certain quadrature rules for Hamiltonian systems with linear operator~$J$, cPG methods are  energy-preserving at the corresponding quadrature nodes~\cite{Gross2005}. 

For a specific class of Hamiltonian and gradient systems in \cite{EggerHabrichShashkov2021} the authors design an energy-consistent Petrov--Galerkin method of arbitrary order with the beneficial property that Galerkin projection in space preserves the Hamiltonian or gradient structure.  
A range of nonlinear systems can be formulated in this framework, and it applies also to some DAEs. 
Unfortunately, for some systems it is not straightforward to see whether they can be reformulated to fit into the framework of \cite{EggerHabrichShashkov2021}, since this involves inverting certain operators. 
For example this is the case for the quasilinear wave equation with friction, see Example~\ref{ex:wave-eq} below. 
A more detailed discussion is presented in Remark~\ref{rmk:comp-EHS} below. 

In~\cite[Sec.~7.4]{morandin24-modeling} a cPG scheme for an alternative formulation of port-Hamiltonian systems is used, which is energy-consistent and allows for a structure preserving space discretization, as presented in \cite[Sec.~7.2]{morandin24-modeling}. 
More specifically, instead of $R(\mathcal{H}'(z))$ as in~\eqref{eq:ph-intro-2},  the author considers terms of the form $\tilde R(z) \mathcal{H}'(z)$, where $\tilde R(z)$ is a $z$-dependent nonlinear operator acting on  $\mathcal{H}'(z)$, and analogously for the term involving $J$.  
While many nonlinear port-Hamiltonian systems fit into both formulations, there are examples of models for which our formulation is more natural since it clearly shows the monotonicity of the operators.  
Examples are the $p$-Laplace equation or the quasilinear wave equation with nonlinear friction term. 
Our strategy to achieve energy consistency and structure preservation by means of a suitable projection of $\mathcal{H'}(z)$ is similar to the one in \cite{morandin24-modeling}, however, our setup avoids problems with singularities in cases where $\tilde R(z)$ is singular, but $\tilde R(z) \mathcal{H}'(z)$ is not, see Remark \ref{rmk:discr-comp} for details.

 There are more classes of methods that are exactly energy-consistent only for the special case of quadratic Hamiltonian, such as standard collocation methods and certain types of Runge--Kutta methods. 
They have been applied to port-Hamiltonian systems including port-Hamiltonian DAEs, which creates extra challenges. 
For example collocation methods were used in  \cite[Thm.~2]{kotyczka19-discrete}, \cite[Sec.~III C]{mehrmann19-structure}, and some
Runge--Kutta methods were used in~\cite{morandin24-modeling}.  

To summarize, before the first announcement of this work the only high-order energy-consistent and structure-preserving methods available for (port)-Hamiltonian systems with non-quadratic Hamiltonian are energy-preserving collocation methods~\cite{hairer10-energy} for general port-Hamiltonian systems, and continuous Petrov--Galerkin methods for a specific class of Hamiltonian and gradient systems~\cite{EggerHabrichShashkov2021}, and \cite[Sec.~7]{morandin24-modeling} projecting in a different way, see also Remark~\ref{rmk:classical} below for a comparison of the formulations.  

In this work we develop a (high-order) energy-consistent cPG method for port-Hamiltonian systems  of the form~\eqref{eq:ph-intro-2}, including Hamiltonian systems for general Hamiltonian.  
More specifically, we present a modified continuous Petrov--Galerkin method that coincides with the classical one for the special case of quadratic Hamiltonian and linear $J$ and $R$. 
However, it is energy-consistent also for non-quadratic Hamiltonian and nonlinear $J$ and $R$.  
Even though the type of discretization is similar to the one presented in~\cite[Sec.~7]{morandin24-modeling} they agree only if $\tilde J(z), \tilde R(z)$ are independent of $z$ and  $J,R$ are linear, but not in general. 
See also Remark~\ref{rmk:discr-comp} below for a detailed discussion. 
Our framework is not suitable for problems where irregular solutions display additional dissipation effects, e.g., shock solutions in hyperbolic models. 
Also, our approach does not apply to models with degeneracy which induces lack of regularity. 
In combination with a suitable space discretization the (port)-Hamiltonian structure is preserved under space discretization. 
 
In a very recently announced and independent work~\cite{AndrewsFarrell2024} the authors construct a general class of high-order in time energy-preserving methods for Hamiltonian systems and dissipative differential equations. 
Remarkably, they even achieve the conservation of an arbitrary number of invariants by use of an additional auxiliary variable per invariant. 
This leads to a larger mixed problem and can be reformulated in terms of $L^2$-projections. 
This approach is showcased on several examples including the compressible Navier--Stokes equation. 
Without control and for state-independent $J$ and $R$ their approach applied to the continuous Petrov--Galerkin agrees with our proposed method. 

\subsubsection*{Main contributions and outline}
The class of nonlinear infinite-dimensional port-Hamiltonian systems under consideration is introduced in Section~\ref{sec:prob-general} in Assumption~\ref{as:problem-general}. 
Examples fitting into this framework include the quasilinear wave equation, possibly with friction and viscosity, and doubly nonlinear parabolic equations, that in special cases reduce to the porous medium equation and the $p$-Laplace equation. 
Those examples are discussed in Section~\ref{sec:ex}. 
Furthermore, in Section~\ref{sec:prob-findim} we present the corresponding finite-dimensional setting which encompasses many classical finite-dimensional port-Hamiltonian systems, and we present some examples in Section~\ref{sec:ex}. 
However, our framework does not include systems with state-dependent operator $C$ and DAEs. 

For the class of port-Hamiltonian systems described above, in Section~\ref{sec:scheme} we introduce a continuous Petrov--Galerkin scheme (Scheme~\ref{def:scheme}) of arbitrary polynomial degree~$k \in \mathbb{N}$.
On the nonlinear terms a quadrature formula is used, and we only require positivity of the quadrature weights. 
Our scheme is designed to be energy-consistent, as proved in Proposition~\ref{prop:discr-energy}. 
The main tool to achieve this for general Hamiltonian is the $L^2$-projection mapping to piecewise polynomials of maximal degree~$k-1$ in time. 
This is rather natural, since stability estimates for the $L^2$-projection of the solution can be obtained as in~\cite[Lem.~4.1]{Ahmed2015}. 
For linear port-Hamiltonian systems with quadratic Hamiltonian the projection cancels, and hence our proposed scheme reduces to the standard cPG method. 

Our approach has similarities with mixed methods, because it can be reformulated by use of an auxiliary variable. 
It has been extended in the more recent work~\cite{AndrewsFarrell2024}. 
For finite-dimensional port-Hamiltonian systems similar ideas were used in \cite{chaturantabut16-structure} for the purpose of structure preserving model reduction. 
For infinite-dimensional (port)-Hamiltonian systems the use of an orthogonal projection in space allows to preserve the (port)-Hamiltonian structure under space discretization. 
In the PhD thesis~\cite{Jackaman2019} such an approach was used for a DG space discretization of some dispersive Hamiltonian system, and was extended in several contributions thereafter~\cite{CelledoniJackaman2021,JackamanPryer2021} for continuous and discontinuous space discretization for several Hamiltonian PDEs. 
In Section~\ref{sec:space-discr} we present the corresponding structure preserving conforming space discretization in our general framework. 
A similar discretization has also been used for port-Hamiltonian systems of a slightly different form and for port-Hamiltonian DAEs in~\cite[Sec.~7.2]{morandin24-modeling}.  

Furthermore, in Section~\ref{sec:num-exp} we present numerical experiments using Gauß quadrature to verify the energy consistency for some of the examples presented in Section~\ref{sec:ex}. 
Indeed, we need quadrature to handle the projections as well as the nonlinear terms involving~$R$, $J$ and $B$.
More specifically, we examine for which type of quadrature the energy consistency is satisfied up to machine precision. 
Additionally, we investigate the convergence in the time discretization parameter. 
We observe that using Gauß quadrature with $k$ nodes (with exactness degree $2k-1$) both for the terms involving $J$ and $R$ and for the approximation of the projection, leads to optimal convergence rate $k+1$ in the time discretization parameter of the error in the $L^\infty$-norm. 
Of course, this requires sufficient regularity of the solutions. 
This matches with the results in~\cite{Aziz1989} and confirms that the use of the $L^2$-projection does not corrupt the convergence if applied correctly. 
Furthermore, we observe superconvergence at the time grid points of order $2k$ in the time discretization for some of the examples. 
These are exactly the rates known for parabolic problems~\cite{Aziz1989}, and hence also in this regard we do not lose anything by means of the projection.  
Moreover, for some infinite-dimensional examples we use a structure preserving space discretization. 
We find that the projection does not reduce the convergence order. 
Additionally, we investigate the convergence order in the time discretization parameter for several space discretization parameters, and observe that the errors are independent of the space discretization. 
Energy consistency is achieved by using a sufficiently high number of nodes in the Gauß quadrature for the projection.
  
\subsubsection*{Discussion and outlook}
The main novelty of our scheme is, that it yields energy consistency for arbitrary Hamiltonian, it preserves the port-Hamiltonian structure under space discretization for a general class of nonlinear infinite-dimensional port-Hamiltonian systems, and it is of arbitrary order.   
It applies to a large class of (port)-Hamiltonian systems of form~\eqref{eq:ph-intro-2}. 
In contrast to~\cite{EggerHabrichShashkov2021} it does not require a reformulation of the model, and in contrast to \cite{morandin24-modeling} it allows to treat port-Hamiltonian systems with singularities. 

Our method extends the standard cPG method in a very natural way. 
Unlike the method in~\cite{EggerHabrichShashkov2021}, our scheme does not coincide with any of the collocation methods~\cite{cohen11-linear,	celledoni17-energy} for special port-Hamiltonian systems, except for polynomial degree $k = 1$. 

Besides proving energy consistency, in this work we do not perform any numerical analysis such as proving well-posedness of the discrete solutions. 
This is due to the fact that nonlinear problems require a highly problem adapted approach, which does not fit well with the general framework we present here. 
Also for a~priori and a~posteriori error estimates more structure of the specific problems has to be used, and is therefore not addressed here. 

Furthermore, port-Hamiltonian systems with non-trivial operator $C(z)$, as well as DAEs and their energy-consistent approximation with cPG methods is not addressed in this work. 
Note that the energy-consistent methods in~\cite{EggerHabrichShashkov2021,morandin24-modeling} are applied also to some DAEs.
See also~\cite{Altmann2022} for a (non energy-consistent) cPG method for certain DAEs. 
To investigate under which structural assumptions or with which modifications the scheme may be applicable to those problems as well, is of great interest, but this is left to future work.

\subsubsection*{Notation.}
Let $X$ and $Y$ be Banach spaces and $I\subseteq \mathbb{R}$ be an interval.
We denote the set of all bounded linear operators mapping from $X$ to $Y$ by $L(X,Y)$. 
The Fréchet derivative of a mapping $f\colon X \to Y$ is denoted by $f'$. 
For an open set $D \subset \BS$ we denote by $C^1(D;\setR)$ the set of mappings, that are Fréchet differentiable at all $v \in D$, and for which 
$D \ni v \mapsto f'(v) \in L(X,Y)$ is continuous. 
We denote the dual space of $X$ by $X'$ and the dual pairing between $X$ and $X'$ by $\langle \cdot, \cdot \rangle_{X',X}$. 

For a bounded Lipschitz domain $\Omega \subset \setR^d$ with $d\in \mathbb{N}$, and $p \in [1,\infty]$ we denote by $L^p(\Omega)$ the Lebesgue space and by $W^{1,p}(\Omega)$ the Sobolev space. 
By $L^p(I;X)$ we denote the Bochner space of functions $g \colon I \to X$ that are Bochner measurable and that have finite norm, with norms
\begin{align*}
  \norm{g}_{L^p(I;X)} &\coloneqq \left(\int_I \norm{g(t)}_X^p \dt\right)^{1/p} \text{ if } p \in [1, \infty),\\
  \norm{g}_{L^\infty(I;X)} &\coloneqq \esssup_{t \in I}\norm{g(t)}. 
\end{align*}
Furthermore, $C(I;\BS)$ denotes the space of continuous functions with values in $\BS$. 
Generic constants $c$ may change in a sequence of inequalities and only depend on the quantities specified. 


\section{Problem setting}\label{sec:ph-systems}
\noindent 
In this section we introduce the precise formulation of the nonlinear port-Hamiltonian systems under consideration. 
We start with the infinite-dimensional case in Section~\ref{sec:prob-general}, and we present the special case of finite-dimensional port-Hamiltonian systems in Section~\ref{sec:prob-findim}. 
The power balance for these systems is discussed in Section~\ref{sec:pow-bal}. 
Finally, in order to demonstrate the strength of our framework in Section~\ref{sec:ex} we present a number of examples that fit in. 

\subsection{Infinite-dimensional case.}\label{sec:prob-general}
\noindent In the following we consider a general setup for port-Hamiltonian systems. 
For this purpose let $Z$ be a Hilbert space, that we identify with its dual space $Z'$ with inner product $\skp{\cdot}{\cdot}$. 
Let $\BS \hookrightarrow Z$ be a Banach space with dual space $\BS'$ and duality relation $\skp{\cdot}{\cdot}_{\BS',\BS}$. 
Furthermore, for given $T>0$ let $I = [0,T] \subset \setR$ be a time interval. 

We consider port-Hamiltonian systems which in weak form can be stated as 
\begin{subequations}\label{eq:pH-weak}
\begin{align}
	\skp{\partial_t z}{\phi}_{\BS',\BS}  
	&= \bilJ(\eta(z),\phi) 
	- \bilR (\eta(z),\phi)
	 + \bilu(\cdot,\eta(z),\phi) \quad \text{ for all } \phi \in \BS, \text{ pointwise on } I,\\
	z(0) &= z_0, 
\end{align}
\end{subequations}
where $z \colon I \to D \subset \BS$ is a function with $\partial_t z \colon I \to \BS'$, and $z_0 \in \BS$ is a given initial datum. 
We usually shall suppress the time dependence and write for example $z$ instead of $z(t)$. 
The operator $\eta$ is the Fréchet derivative of the Hamiltonian $\mathcal{H}$ and the operators $\bilJ$, $\bilR$ and $\bilu$ are to be specified. 
Here, $\bilJ$ models energy conservative  processes, $\bilR$ describes all dissipative processes, and $\bilu$ contains the control of the system.
For system~\eqref{eq:pH-weak} we make the following assumptions. 

\begin{assumption}[problem formulation]
	\label{as:problem-general}
  \noindent Let $Z$ be a Hilbert space identified with its dual $Z \cong Z'$, with inner product $\skp{\cdot}{\cdot}$, and assume that $(\BS,\norm{\cdot}_{\BS})$ is a reflexive Banach space with continuous and dense embedding $\BS \hookrightarrow Z$. 
  The duality relation between $\BS$ and its dual space $\BS'$ is denoted by $\skp{\cdot}{\cdot}_{\BS',\BS}$. 
	Let $I = [0,T] \subset \setR$ be a bounded interval. 
\begin{enumerate}[label=  (A\arabic*)]
\item  \label{itm:H}
$\mathcal{H} \colon D \to \setR$  is a Fréchet differentiable operator with $\mathcal{H} \in C^1(D;\setR)$ for some open subset $D \subset Y \hookrightarrow \BS$ for a Banach space $Y$ and we denote $\eta \coloneqq \mathcal{H}'$; 
\item \label{itm:forms} 
$\bilJ, \bilR \colon \BS \times \BS \to \setR$ and $ \bilu \colon I \times \BS \times \BS \to \setR$ are functionals, which are linear and Lipschitz continuous in their last argument, and continuous in the remaining arguments. 
We assume that there exists $p\in (1,\infty)$, a function $\mathfrak{b} \in C(I)$ and a constant $c>0$ such that 
\begin{align}\label{est:bd}
		\abs{\bilJ(v,w)} + \abs{\bilR(v,w)} + \abs{\bilu(t,v,w)}
	\leq c( \mathfrak{b}(t) + \norm{v}_{\BS}^{p-1}) \norm{w}_{\BS},
\end{align}
for any $v,w \in \BS$ and any $t \in I$. 
 Furthermore, for $\bilJ,\bilR$ we assume that 
\begin{enumerate}[label = (A2\roman*),leftmargin=1.3cm]
\item \label{itm:bJ}
$\bilJ$ satisfies for any  $v \in \BS$ that  
	\begin{align*}
		\bilJ(v,v) = 0;
	\end{align*}
\item \label{itm:bR-dissip}
  $\bilR$ is dissipative in the sense that for any  $v \in \BS$ we have that 
\begin{align*}
	\bilR(v,v) \geq 0. 
\end{align*}
\end{enumerate}
\end{enumerate}
\end{assumption}

\begin{remark}\hfill
	\begin{enumerate}[label = (\alph*)]\label{rmk:assump-gen}
		\item 
		There are several nonlinearities in~\eqref{eq:pH-weak}: $\mathcal{H}$ need not be quadratic, and hence $\eta$ may be nonlinear, and also $j$, $r$ and $b$ may be nonlinear. 
		\item 
		 DAEs are not contained in the framework, because in~\eqref{eq:pH-weak} there is no (degenerate) operator acting on $\partial_t z$.
	\item  
Notably,~\ref{itm:bJ} is satisfied, if~$\bilJ$ is skew-symmetric.  
\item The estimate \eqref{est:bd} in~\ref{itm:forms} on the functionals ensures that all terms on the right-hand side of~\eqref{eq:pH-weak} are integrable in time for test functions in $L^p(I;\BS)$ and provided that the solution is sufficiently regular. 
 Later we shall consider solutions $z$ to~\eqref{eq:pH-weak} with $\eta(z) \in C(I;\BS)$.  
	\end{enumerate}
\end{remark}

\begin{remark}\label{rmk:comp-EHS}
	In \cite{EggerHabrichShashkov2021} the authors consider Hamiltonian and gradient systems of the form 
	\begin{align}\label{eq:pH-EHS}
		\mathcal{C}(z) \partial_t z = - \mathcal{H}'(z) + f(z),
	\end{align}
	for Hamiltonian $\mathcal{H}$, and $\mathcal{C}(\cdot)$ positive semi-definite for any argument, and for given function $f$. 
	This formulation has the advantage, that the structure is preserved under Galerkin projection in space, which is not the case for systems of the form~\eqref{eq:pH-weak}, as discussed in detail in~\cite{EggerHabrichShashkov2021}. 
	For invertible $\mathcal{C}(z)$, one may transform the systems~\eqref{eq:pH-weak} and~\eqref{eq:pH-EHS} into each other if either for given $\mathcal{C}$ the operators $j,r$ are such that 
	\begin{align*} 
		- \skp{\mathcal{C}(z)^{-1} \eta(z)}{\phi} = \bilJ(\eta(z),\phi) - \bilR(\eta(z),\phi),
	\end{align*}
	or, if for given $r,j$, the operator $\mathcal{C}$ is such that the same identity holds.  
	Below in Section~\ref{sec:ex} we present examples of infinite-dimensional systems that fit in our framework.   
	However, for one of them (Example~\ref{ex:wave-eq} on the wave equation with viscosity or friction) a reformulation in the form~\eqref{eq:pH-EHS} is not available to the best of our knowledge. 
	Indeed,  the difficulty consists in the fact that the right-hand side has several terms that would require a different operator $\mathcal{C}^{-1}$ each. 
	Note that for two invertible operators their sum is not invertible in general. 
	Even for linear operators on finite-dimensional spaces this  requires certain commutation properties. 
\end{remark}

\subsection{Finite-dimensional case}\label{sec:prob-findim}
\noindent 
For the purpose of visualization let us also present the port-Hamiltonian system for the simpler case with finite-dimensional (Hilbert) space $Z = \BS = \setR^d$ for some $d \in \mathbb{N}$ with the Euclidean inner product $\skp{\cdot}{\cdot}_{\ell^2}$. 
Then, the port-Hamiltonian system~\eqref{eq:pH-weak} reduces to
\begin{subequations}\label{eq:pH-findim}
	\begin{align}
 \partial_t z
		&= J(\eta(z))- R(
		\eta(z))
		+ B(\cdot,\eta(z)) \quad \text{ pointwise on } I,\\
		z(0) &= z_0,
	\end{align}
\end{subequations}
for $z_0 \in \setR^d$. 
In this special case Assumption~\ref{as:problem-general} simplifies considerably. 

\begin{assumption}[problem formulation in finite dimensions]
	\label{as:problem-findim} 
\noindent With $d \in \mathbb{N}$ we assume that: 
\begin{enumerate}[label=  (a\arabic*)]
		\item  \label{itm:H-v}   
	$\mathcal{H} \in C^1(D;\setR)$ for some open subset $D \subset \setR^d$ and we denote $\eta \coloneqq \mathcal{H}'$; 
	\item \label{itm:forms-v}
	$J,R\colon \setR^d \to \setR^d$ and $B \colon I\times \setR^d  \to \setR^d$ are continuous mappings. 
	We assume that there exists $p\in (1,\infty)$, a function $\mathfrak{b} \in C(I)$ and a constant $c>0$ such that 
	\begin{align}\label{est:bd-v}
		\norm{J(v)}_{\ell^2} + \norm{R(v)}_{\ell^2} +  \norm{B(t,v)}_{\ell^2}
		\leq c( \mathfrak{b}(t) + \norm{v}_{\ell^2}^{p-1}),
	\end{align}
	for any $v\in \setR^d$ and any $t \in I$. 
	 Furthermore, on $J$ and $R$ we assume that
	\begin{enumerate}[label = (a2\roman*),leftmargin=1.3cm]
		\item \label{itm:J}
		$J$ satisfies for any $v \in \setR^d$ that 
		\begin{align*}
			\skp{J(v)}{v} = 0;
		\end{align*}
		\item \label{itm:R}
		$R$ is dissipative in the sense that for any $v \in \setR^d$ one has 
		\begin{align*}	
			\skp{R(v) }{v}_{\ell^2} \geq 0.
		\end{align*}
	\end{enumerate}	
\end{enumerate}
\end{assumption}

\begin{remark}[classical state dependent port-Hamiltonian systems]\label{rmk:classical}
	\hfill 
	The rather general conditions on $J, R, B$ in~\ref{itm:forms-v} include $z$-dependent operators, if $\eta$ is invertible. 
	Note however, that in~\eqref{eq:pH-findim} there is no $z$-dependent operator acting on $\partial_t z$, and hence no DAEs are contained in the framework. 
	 	
	In this situation the classical formulation of port-Hamiltonian systems reads
	\begin{subequations}\label{eq:pH-findim-classical}
		\begin{align}
	 \partial_t z
			&= \left(\tilde J(z)- \tilde R(z)\right)\eta(z)
			+ \tilde B(z) u \quad \text{ pointwise on } I,\\
			z(0) &= z_0. 
		\end{align}
	\end{subequations}
	Here~$u\colon I \to \setR^l$ with~$l \leq d$ is some control and~$z_0 \in \setR^d$ is a given initial datum.
	For each $v\in \setR^d$ one has that $\tilde J(v), \tilde R(v) \in \setR^{d \times d}$ and that $\tilde B(v) \in \setR^{d \times l}$. 

	If $\eta \colon \setR^d \to \setR^d$ is invertible, then $J, R \colon \setR^d \to \setR^d$ and $B \colon I \times \setR^d \to \setR^d$ can be defined by 
	\begin{align*}
	J(v) \coloneqq \tilde J(\eta^{-1}(v))v, \quad 
	R(v) \coloneqq \tilde R(\eta^{-1}(v))v 
	 \quad  \text{ and } \quad 
	B(t,v) \coloneqq \widetilde{B}(\eta^{-1}(v)) u(t)
	\end{align*}	
	 for any $v \in \setR^d$, and any $t \in I$. 
	Indeed, for Assumption~\ref{as:problem-findim} to be satisfied it suffices that the following conditions hold: 
	\begin{enumerate}[label = (a\arabic*')] 
	\item \label{itm:H-c}
	$\mathcal{H}  \in C^2(\setR^d;\setR)$ is  strictly convex and we denote $\eta \coloneqq \mathcal{H'} \in C^1(\setR^d;\setR^d)$. 
	\item \label{itm:forms-c}
	$\tilde J(v)$, $\tilde R(v) \in \setR^{d \times d}$ and $\tilde B(v) \in \setR^{d\times l}$ are continuous as functions in $v \in \setR^d$. 	
	We assume that there exists $p\in (1,\infty)$, 
	 and a constant $c>0$ such that 
			\begin{align}\label{est:bd-c}
			\norm{\tilde J(v)}_{2} + \norm{\tilde R(v)}_{2}+ \norm{\tilde B(v)}_{2} 
			\leq c( 1 + \norm{\eta(v)}_{\ell^2}^{p-2}),
		\end{align}
		for any $v \in \setR^d$. 
	Here $\norm{\cdot}_2$ denotes the matrix norm induced by $\norm{\cdot}_{\ell^2}$. 
	 Furthermore, we assume that
	\begin{enumerate}[label = (a2\roman*'),leftmargin=1.3cm]
			\item \label{itm:J-c}
			for any $v \in \setR^d$ the matrix $\tilde J(v)$ is skew-symmetric;
			\item \label{itm:R-c}	
			for any $v \in \setR^d$ the matrix $\tilde R(v)$ is positive semi-definite. 
		\end{enumerate}
	\end{enumerate}
The condition that $\mathcal{H}$ is strictly convex in~\ref{itm:H-c} indeed implies, that $\eta' = \mathcal{H}'' >0$, and hence $\eta$ is invertible. 
	
Note that in the situation of~\eqref{eq:pH-findim-classical} the splitting of an operator into a (symmetric) positive definite part and a skew-symmetric part is unique. 
In the more general case of genuinely nonlinear operators in finite dimensions~\eqref{eq:pH-findim}, and  in infinite dimensions~\eqref{eq:pH-weak} a splitting into dissipative part $r$ and a conservative part $j$ is not so obvious, and depends on the modeling. 
However, a unique splitting of the corresponding right-hand side of~\eqref{eq:pH-weak} is not used in the following investigation. 
	
In~\cite{morandin24-modeling} finite and infinite-dimensional systems of the form~\eqref{eq:pH-findim-classical} as well as the corresponding DAEs are investigated. 
Even though the examples considered below fit in both formulations, some of the resulting energy-consistent discretizations differ, as explained in Remark~\ref{rmk:discr-comp} below. 
\end{remark}

\begin{remark}[structure preserving space discretization]
	Our approach to achieve an energy-consistent time discretization of port-Hamiltonian systems of the form~\eqref{eq:pH-weak}, and \eqref{eq:pH-findim} will be presented in Section~\ref{sec:scheme}. 
	For infinite-dimensional systems it is closely related to a structure-preserving space discretization, see Section~\ref{sec:space-discr} below. 
	With this approach the discrete system is a finite dimensional port-Hamiltonian system of form~\eqref{eq:pH-findim}. 
\end{remark}

\subsection{Power balance}\label{sec:pow-bal}
\noindent The structural assumptions on the port-Hamiltonian system  entail that the Hamiltonian is diminished by the dissipative term $\bilR$, and affected by $\bilu$ which may include boundary terms and controls. 

\begin{lemma}[power balance]\label{lem:pow-bal}
	Let Assumption~\ref{as:problem-findim} be satisfied. 
	Then, any sufficiently smooth solution $z$ to  system~\eqref{eq:pH-weak}, with $\eta(z) \in C(I;\BS)$, satisfies the following \emph{power balance}
	\begin{equation}\label{eq:powerbalance}
		\begin{aligned}
		\frac{\mathrm{d}}{\mathrm{d}t} \mathcal{H}(z) 
		&= - \bilR(\eta(z),\eta(z)) + \bilu(\cdot,\eta(z),\eta(z))&& \\
		& \leq  \bilu(\cdot,\eta(z),\eta(z)) \qquad &&\text{ pointwise on } I.
		\end{aligned}
	\end{equation}
\end{lemma}
\begin{proof}
	Using the properties of~$\cH$  in Assumption~\ref{as:problem-general}~\ref{itm:H}, employing equation~\eqref{eq:pH-weak} and condition~\ref{itm:bJ} on $j$ we find that 
	\begin{align*}
		\frac{\mathrm{d}}{\mathrm{d}t} \cH(z) 
		& = \skp{\partial_t z}{\cH'(z)}_{\BS',\BS}
		= \skp{ \del_t z}{ \eta(z)}_{\BS',\BS}\\
		& = \bilJ(\eta(z),\eta(z)) - \bilR(\eta(z),\eta(z)) + \bilu(\cdot,\eta(z),\eta(z)) \\
		& =  - \bilR(\eta(z),\eta(z)) + \bilu(\cdot,\eta(z),\eta(z)). 
	\end{align*}
In combination with the dissipative nature of $\bilR$ according to~\ref{itm:bR-dissip} this proves the claim. 
\end{proof}

\begin{remark}\hfill
\begin{enumerate}[label = (\alph*)]
	\item	An alternative formulation of~\eqref{eq:powerbalance} is the \emph{energy balance}
	\begin{equation}
	\label{eq:energybalance}
	\begin{aligned}
		\cH(z(t_1)) - \cH(z(t_0)) &= \int_{t_0}^{t_1} - \bilR(\eta(z),\eta(z)) + \bilu(\cdot,\eta(z),\eta(z))\dt \quad  \\
		 &\leq \int_{t_0}^{t_1}  \bilu(\cdot,\eta(z),\eta(z))\dt \qquad \text{ for any }t_0, t_1 \in I.
		 \end{aligned}
	\end{equation}
\item For finite-dimensional port-Hamiltonian systems~\eqref{eq:pH-findim} the power-balance analogously reads 
\begin{align}\label{eq:powerbal-1}
	\frac{\mathrm{d}}{\mathrm{d}t} \mathcal{H}(z) 
= - \skp{R(\eta(z))}{\eta(z)}_{\ell^2}+ \skp{B(\cdot,\eta(z))}{\eta(z)}_{\ell^2}
\leq   \skp{B(\cdot,\eta(z))}{\eta(z)}_{\ell^2}. 
\end{align}
For a finite-dimensional port-Hamiltonian system of form~\eqref{eq:pH-findim-classical} the energy balance reads
\begin{align}\label{eq:powerbal-2}
	\frac{\mathrm{d}}{\mathrm{d}t} \mathcal{H}(z) 
	= - \skp{\tilde R(z) \eta(z)}{\eta(z)}_{\ell^2}+ \skp{\tilde B(z)u}{\eta(z)}_{\ell^2}
	\leq  \skp{\tilde B(z)u}{\eta(z)}_{\ell^2}. 
\end{align}
In the context of port-Hamiltonian modeling the quantity $y \coloneqq \tilde B(z)\tp \eta(z)$ is the system output. 
In this case $\skp{\tilde B(z) u}{\eta(z)}_{\ell^2}$ can be replaced by $\skp{y}{u}_{\ell^2} $ to arrive at the customary formulation of the energy  balance for port-Hamiltonian systems.

Note that in the finite-dimensional case the condition $\eta(z) \in C(I;\BS)$ is satisfied, provided that $z \in C(I;\setR^d)$. 
	\end{enumerate}
\end{remark} 
 
 \subsection{Examples}\label{sec:ex}
\noindent Let us discuss examples of port-Hamiltonian systems of the form~\eqref{eq:pH-weak} 
that fit into the framework presented in Assumption~\ref{as:problem-general}. 
We focus on cases, where the system is nonlinear or the Hamiltonian is not quadratic. 
We start by considering some finite-dimensional examples, that fit into the framework in Section~\ref{sec:prob-findim}. 
Then, we proceed with some infinite-dimensional examples for systems as introduced in Section~\ref{sec:prob-general}. 

Since the purpose of this section is to showcase the strength of the framework we usually refrain from presenting optimal estimates. 
In particular, the estimates on $r,j,b$ as in Assumption~\ref{as:problem-general}~\ref{itm:forms} can be sharpened in several places. 
Some pointers of generalizations are given in Remark~\ref{rmk:power-extensions}.  
\medskip 

First, we consider a finite-dimensional system for which $\tilde J$ and $\tilde R$ are independent of $z$, but the Hamiltonian is non-quadratic.  

\begin{example}[Toda lattice]\label{ex:toda}
	The Toda lattice describes the motion of a chain of particles in 1D, where each particle is connected to its nearest neighbors with an exponential spring,  cf.~\cite[Sec.~3.4.2]{chaturantabut16-structure}.
	For $N \in \mathbb{N}$ the number of particles, $q\in \setR^N$ the displacement vector of the particles and $p \in \setR^N$ the momentum vector of the particles we set  
	\begin{align*}
		z \coloneqq \begin{pmatrix} 
			q \\ p \end{pmatrix} \in \setR^{2N}.
	\end{align*}
	The following system of ordinary differential equations describes the motion 
	\begin{equation}\label{eq:toda}
		\del_t z = (\tilde J - \tilde R)\eta(z) +  \tilde Bu, \quad \
	\end{equation}
	for given $u \in \setR$, 
	 and for matrices $\tilde J$, $\tilde R$ and $\tilde B$ given by
	\begin{align*}
		\tilde J = \begin{pmatrix} 0 & \identity_N \\ -\identity_N & 0 \end{pmatrix}  \in \mathbb{R}^{2N \times 2N},
		\quad 
	\tilde	R &= \begin{pmatrix} 0 & 0 \\ 0 & \mathrm{diag}(\gamma_1,\dots,\gamma_N) \end{pmatrix} \in \mathbb{R}^{2N \times 2N}, \\\text{ and } 
	\tilde	B &= \begin{pmatrix}  0 \\  e_1 \end{pmatrix} \in \mathbb{R}^{2N \times 1}.
	\end{align*}
	Here $\identity_N \in \setR^{N \times N}$ is the identity matrix, $\gamma_i\geq 0$ are given damping parameters, and $e_1$ is the first unit vector. 
	This means that the control exerts a force on the first particle. 
	The Hamiltonian of the system reads 
	\begin{equation}\label{def:H-toda}
		\cH(z) = \sum_{k=1}^{N} \frac12 p_k^2 + \sum_{k=1}^{N-1} \exp(q_k - q_{k+1}) + \exp(q_N - q_1) - N \quad \text{ for } 	z \coloneqq \begin{pmatrix} 
			q \\ p \end{pmatrix}, 
	\end{equation}
	which is smooth and strictly convex and $\eta \coloneqq \mathcal{H}'$. 
	Since $\tilde J$, $\tilde R$ and $\tilde B$ are independent of $z$, the finite-dimensional system~\eqref{eq:toda} is of the form~\eqref{eq:pH-findim-classical} and it suffices to verify~\ref{itm:H-c}--\ref{itm:forms-c} in Remark~\ref{rmk:classical}. 
	One can check directly that~\ref{itm:H-c} holds. 
	The matrices $\tilde J, \tilde R, \tilde B$ are independent of $z$, and hence continuous in $z$. 
	The estimate~\eqref{est:bd-c} holds with $p = 2$. 
	Since $\tilde J$ is skew-symmetric and $\tilde R$ is positive semi-definite,~\ref{itm:forms-c} is satisfied. 
	By standard ODE theory well-posedness is available. 
\end{example}

Next, let us consider a finite-dimensional system with quadratic Hamiltonian, but $z$-dependent operator $\tilde J $. 

\begin{example}[spinning rigid body]\label{ex:spinning-body}
	In~\cite[Example~6.2.1]{vanderschaft17-l2gain} a rigid body spinning around its center of mass in the absence of gravity is modeled by a port-Hamiltonian system of the form
	\begin{equation}\label{eq:rigid-body}
		\del_t z = \tilde J(z)Qz + \tilde Bu. 
	\end{equation}
	Here, the solution $z = (p_1, p_2, p_3)\tp \in \mathbb{R}^{3}$ is the vector of the angular momenta of the body in the three spatial dimensions. 
	The matrices $\tilde J(z), Q \in \setR^{3 \times 3}$ and $\tilde B \in \setR^{3 \times 1}$ are given by 
	\begin{equation*}
		\tilde J(z) = \begin{pmatrix} 0 & -p_3 & p_2 \\ p_3 & 0 & -p_1 \\ -p_2 & p_1 & 0 \end{pmatrix}, ~
		Q = \begin{pmatrix} \frac1{I_1} & 0 & 0 \\ 0 & \frac1{I_2} & 0 \\ 0 & 0 & \frac1{I_3} \end{pmatrix},\quad  \text{ and } \quad 
	\tilde	B = \begin{pmatrix} b_1 \\ b_2 \\ b_3 \end{pmatrix},
	\end{equation*}
	with $I_1,I_2,I_3>0$ the principal moments of inertia, $\tilde B$ contains the coordinates $b_1, b_2, b_3 \in \setR$ of the axis around which torque is applied, and $u\in \setR$ is a given control.
	The system~\eqref{eq:rigid-body} is of the form~\eqref{eq:pH-findim-classical} with $\tilde R \equiv 0$, with quadratic Hamiltonian $\mathcal{H}(z) = \tfrac{1}{2} z\tp Q z$ and $\eta(z) \coloneqq  \mathcal{H}'(z) = Qz$. 
	Again, to verify Assumption~\ref{as:problem-findim} it suffices to check~\ref{itm:H-c}-\ref{itm:forms-c}. 
	Indeed, condition~\ref{itm:H-c} is immediate. 
	Furthermore, $\tilde J(z)$ is skew-symmetric, linear in $z$, and in particular continuous.
	The matrix $\tilde B$ is independent of $z$ and hence also continuous.  
	Thanks to
	\begin{align*}
		\norm{\tilde J(z)}_2 + \norm{\tilde B}_2 \leq c (1 +  \norm{z}_{\ell^2}) \leq c (1 + \norm{\eta(z)}_{\ell^2})\quad \text{ for any } z \in \setR^3,
	\end{align*}
	with constant depending on the dimension $3$, on $\tilde B$ and on $Q$,~\ref{itm:forms-c} is satisfied with $p = 3$. 
	Note that $u$ is independent of time. 	
By standard ODE theory well-posedness of solutions is available.
\end{example}
 
 Let us now proceed to consider examples of infinite-dimensional port-Hamiltonian systems. 
 They all are evolution equations posed on~$Q = I \times \Omega$ for given final time~$T > 0$, time interval~$I \coloneqq [0,T]$ and for a bounded Lipschitz domain $\Omega \subset \setR^d$, with $d\in \mathbb{N}$. 

\begin{example}[quasilinear wave equation]\label{ex:wave-eq}
We consider the quasilinear wave equation with friction and viscosity: 
	\begin{equation}\label{eq:wave-general}
		\begin{aligned}
			\partial_t \rho + \div(v) & = 0, \\
			\partial_t v + \div (p(\rho)) & = - \gamma(\rho) F(v) +\nu \Delta v,
		\end{aligned}  
	\end{equation}
	on $Q$, 
	subject to initial conditions~$(\rho,v)(0,\cdot) =(\rho_0, v_0)$ in~$\Omega$. 
	The system is supplemented by suitable  boundary conditions introduced in the sequel. 
	
	Here,~$p, \gamma$ and $F$ are given, possibly nonlinear functions to be specified, and~$\nu \geq 0$ is a given parameter. 
	The term~$\gamma(\rho) F(v)$ represents friction forces, i.e., there is no friction if this term vanishes. 
	Similarly, the term~$\nu \Delta v $ represents viscous forces, and the non-viscous case is recovered for~$\nu = 0$.  	
	The function $p\colon \setR \to \setR$ 
	is assumed to be strictly monotone and continuously differentiable, 
	 and $\gamma \colon \setR \to [0,\infty)$ is sufficiently smooth and bounded.  
	Typical examples for $F \colon \setR \to \setR$ are power law type functions, and in the following we consider
	\begin{align}\label{eq:friction-powerlaw}
		F(v) = \abs{v}^{s-2} v \quad  \text{ for some } s \in (1,\infty),
	\end{align}	
	see Remark~\ref{rmk:power-extensions} for possible extensions. 

	Let us verify that the weak formulation of~\eqref{eq:wave-general} can be cast into the form~\eqref{eq:pH-weak} and that Assumption~\ref{as:problem-general} holds. 
	We choose $Z = L^2(\Omega)^{1+d}$ and  set  $z_0 = (\rho_0,v_0)\tp$ and 
	\begin{align}\label{eq:wave-z}
		z \coloneqq \begin{pmatrix} \rho \\ v \end{pmatrix}.
	\end{align} 
	For some constant $p_0 \in \setR$ we choose as Hamiltonian 
	\begin{align}\label{eq:H-wave}
		\mathcal{H}(z) \coloneqq \int_{\Omega} P(\rho) \dx +  \frac{1}{2} \int_{\Omega} \abs{v}^2 \dx,
		   \qquad \text{ for } P(\rho) \coloneqq \int_{0}^\rho p(r) \, \mathrm{d}r + p_0, 
	\end{align}
	Note that $p$ is strictly monotone if and only if $\mathcal{H}$ is strictly convex. 
	Then, formally, one has 
	\begin{align}\label{eq:wave-eta}
 \eta(z) \coloneqq \mathcal{H}'(z) = \begin{pmatrix} p(\rho) \\ v \end{pmatrix}.
	\end{align} 
	Multiplying the system~\eqref{eq:wave-general} by a smooth test function $\phi = (\xi , w)\tp \in C^\infty(Q)^{1+d}$, assuming smoothness of the solution $z$,  integrating over $\Omega$, and integration by parts on the second and on the last term on the right-hand side we obtain 
	\begin{align}
	\int_{\Omega} \partial_t z \cdot \phi \dx 
	&=  
	- \int_{\Omega} \div(v) \, \xi\dx 
	 - \int_{\Omega} \nabla (p(\rho)) \cdot w \dx 
	 -  \int_{\Omega} \gamma(\rho) F(v) \cdot w \dx 
	 + \nu \int_{\Omega} \Delta v \cdot w \dx
	 \notag
	 \\
	& =	- \int_{\Omega} \div(v) \, \xi\dx 
	-  \int_{\partial \Omega }p(\rho)\, w  \cdot n \dsigma
	+  \int_{\Omega}  p(\rho)\, \div(w) \dx 
	\label{eq:wave-wform1}
	\\
	& \qquad -  \int_{\Omega} \gamma(\rho) \abs{v}^{s-2}
	v \cdot w \dx + \nu \int_{\partial \Omega} w\tp \nabla v \, n \dsigma 
	- \nu \int_{\Omega} \nabla v : \nabla w \dx, 
		\notag 
	\end{align} 
where $:$ denotes the Frobenius product between two matrices, and $n$ is the outer unit normal to~$\partial \Omega$. 
This weak formulation is to be understood as an identity pointwise in time. 

To consider the boundary conditions and to identify the Banach space $\BS \hookrightarrow Z$, let us distinguish the cases with and without viscosity: 
\begin{enumerate}[label = (\alph*)]
	\item If $\nu \neq 0$, then we may consider homogeneous boundary conditions on $v$ on part of the boundary $\Gamma$ and natural boundary conditions on the relative complement $\partial \Omega \setminus \Gamma$, i.e., 
	\begin{equation}\label{eq:bc-wave-1}
	\begin{aligned}
		v &= 0 \quad &&\text{ on } I \times \Gamma,\\
	(p(\rho) \identity - \nu \nabla v) n & = g && \text{ on } I \times \partial \Omega \setminus \Gamma,
	\end{aligned}
	\end{equation}
	for some given function~$g \in C(I;L^2(  \partial \Omega \setminus \Gamma)^d)$. 
	We assume that $\Gamma$ is sufficiently regular, e.g., it has finitely many components with sufficiently smooth boundary in $\partial \Omega$. 
	Starting from~\eqref{eq:wave-wform1}, imposing the first condition on the function space and the latter one weakly we arrive at
	\begin{equation}	\label{eq:wave-wform2}
		\begin{aligned}
		\int_{\Omega} \partial_t z \cdot \phi \dx 
		& =	- \int_{\Omega} \div(v) \, \xi\dx 
			+  \int_{\Omega}  p(\rho)\, \div(w) \dx 
			\\
	& \quad 	-  \int_{\Omega} \gamma(\rho) \abs{v}^{s-2} v \cdot w \dx
		- \nu \int_{\Omega} \nabla v : \nabla w \dx
	 - \int_{\partial \Omega \setminus \Gamma } g \cdot w   \dsigma,
	\end{aligned} 
	\end{equation}
	for any	$\phi = (\xi , w)\tp \in C^\infty(Q)^{1+d}$ with $w|_{I \times \Gamma} = 0$. 
	This suggests the choice of function spaces with continuous and dense embedding 
	\begin{align}\label{def:VZ-wave-1}
		\BS \coloneqq L^2(\Omega) \times (H^1_{\Gamma}(\Omega)^d \cap L^s(\Omega)^d) \hookrightarrow L^2(\Omega) \times L^2(\Omega)^d \eqqcolon Z,
	\end{align} 
	with $s\geq 2$ as in~\eqref{eq:friction-powerlaw}. 
	Here $H^1_\Gamma(\Omega)^d$ is defined as closure of $\{f \in C^{\infty}(\Omega)^d \colon f|_{\Gamma} = 0\}$ with respect to the norm in $H^1(\Omega)^d$. 
	Depending on the growth properties of $p$ the Banach space $Y \subset X$ and the open set $D \subset Y$ have to be chosen sufficiently small that $\mathcal{H} \colon D \to \setR$ as in~\eqref{eq:H-wave} is well-defined.
	 Fréchet differentiability of $\mathcal{H}$ as well as continuity of $D \to Y'$,   $v\mapsto \mathcal{H}'(v) \eqqcolon \eta(v)$ are available, provided that $p$ is continuously differentiable and $D\subset Y \hookrightarrow X$ is sufficiently small. 
	Under those conditions Assumption~\ref{as:problem-general}~\ref{itm:H} is satisfied.  	

	We define the functionals $j,r\colon \BS \times \BS \to \setR$ and $b\colon I \times \BS \times \BS \to \setR$ by
	\begin{align}\label{eq:forms-wave-j}
		\bilJ(\chi, \phi)  &\coloneqq - \int_{\Omega} \div(\chi_2) \phi_1  \dx + \int_{\Omega} \chi_1 \div(\phi_2)\dx,  \\
		\label{eq:forms-wave-r}
		\bilR(\chi,\phi) & \coloneqq \nu \int_{\Omega} \nabla \chi_2 : \nabla \phi_2 \dx + \int_{\Omega} \gamma(\chi_1) \abs{\chi_2}^{s-2} \chi_2 \cdot \phi_2 \dx,\\
		\label{eq:forms-wave-b}
		\bilu(t,\chi,\phi) & \coloneqq - \int_{\partial \Omega \setminus \Gamma} g(t) \cdot \phi_2 \dsigma,
	\end{align} 	
	for any $\chi, \phi \in \BS$ and any $t \in I$. 
If $s \geq 2$ as in~\eqref{eq:friction-powerlaw}, using Young's inequality with $s-1 \geq 1$, a trace inequality and boundedness of $\gamma$, one can show that 
	\begin{align}\label{est:forms-wave-1}
		\abs{\bilJ(\chi,\phi)} + \abs{\bilR(\chi,\phi)} + \abs{\bilu(t,\chi,\phi)} 
		\leq c(1 + \norm{g(t)}_{L^2(\partial \Omega \setminus \Gamma)} + \norm{\chi}_{\BS}^{s-1}) \norm{\phi}_{\BS}, 
			\end{align}
	for any $\phi, \chi \in \BS$, with $\BS$ as defined in~\eqref{def:VZ-wave-1}. 
	If on the other hand $s \in (1,2)$, then we obtain with the continuous embedding $H^1(\Omega) \hookrightarrow L^s(\Omega)$ and Young's inequality with $\frac{1}{s-1}>1$ that 
		\begin{align}\label{est:forms-wave-2}
		\abs{\bilJ(\chi,\phi)} + \abs{\bilR(\chi,\phi)} + \abs{\bilu(t,\chi,\phi)} 
		\leq c(1 + \norm{g(t)}_{L^2(\partial \Omega \setminus \Gamma)} + \norm{\chi}_{\BS}) \norm{\phi}_{\BS}, 
	\end{align}
	for any $\phi, \chi \in \BS$. 
	Consequently, the estimate in Assumption~\ref{as:problem-general}~\ref{itm:forms} is satisfied with~$p = \max(s,2)$.  
	Note that $j,r,b$ are linear in $\phi$, and hence by estimate~\eqref{est:forms-wave-1} also Lipschitz continuous in $\phi$. 
	The same arguments hold for all terms depending on $\chi$ except of the last term of $\bilR$, which ensures continuity in $\chi$. 
	Since $s>1$ in~\eqref{eq:friction-powerlaw}, also the last term of $\bilR$ is continuous in $\chi_2$, and hence in $\chi$. 
	
	To verify~\ref{itm:forms} is remains to show~\ref{itm:bJ} and~\ref{itm:bR-dissip}. 
	Directly from the definition of $\bilJ$ we see that $\bilJ(v,v) = 0$ for any $v\in \BS$. 
	The dissipative nature of $\bilR$, namely $\bilR(v,v) \geq 0$ follows from the fact that the power-law relation~\eqref{eq:friction-powerlaw} is monotone and $\gamma(\cdot) \geq 0$. 

	\item 
	If $\nu = 0$, then in~\eqref{eq:wave-wform1} all terms involving $\nu$ vanish, and we consider the boundary conditions 	
		\begin{equation}\label{eq:bc-wave-2}
		\begin{aligned}
		v \cdot n &= 0 \quad &&\text{ on } I \times \Gamma,\\
p(\rho) & = g && \text{ on } I \times \partial \Omega \setminus \Gamma,
		\end{aligned}
	\end{equation}
	for some given function $g \in C(I; H^1( \Omega))$. 
	Again starting from~\eqref{eq:wave-wform1} with $\nu = 0$ we impose the first condition on the function space and the second one weakly and obtain
		\begin{equation}	\label{eq:wave-wform3}
		\begin{aligned}
		\int_{\Omega} \partial_t z \cdot \phi \dx 
		& =	- \int_{\Omega} \div(v) \, \xi\dx 
		+  \int_{\Omega}  p(\rho)\, \div(w) \dx 
		-  \int_{\Omega} \gamma(\rho) \abs{v}^{s-2}v \cdot w \dx\\
		& \quad 	-  \int_{\partial \Omega \setminus \Gamma}g \, w  \cdot n \dsigma
	\end{aligned}
	\end{equation} 
	for any	$\phi = (\xi , w)\tp \in C^\infty(Q)^{1+d}$ with $w \cdot n |_{I \times \Gamma } = 0$ . 
 For $s\geq 2$ as in~\eqref{eq:friction-powerlaw} we choose the spaces
	\begin{align}\label{def:VZ-wave-2}
		\BS \coloneqq L^2(\Omega) \times (H_{\Gamma}(\div;\Omega)^d \cap L^s(\Omega)^d) \hookrightarrow L^2(\Omega) \times L^2(\Omega)^d \eqqcolon Z,
	\end{align} 
	where $H_{\Gamma}(\div;\Omega)$ denotes the closure of $\{f \in C^{\infty}(\Omega)^d \colon (f \cdot n)|_{\Gamma} = 0\}$ with respect to the norm $\norm{f}_{H(\div;\Omega)} \coloneqq \norm{f}_{L^2(\Omega)}+ \norm{\div f }_{L^2(\Omega)}$. 
	The validity of~\ref{itm:H} follows analogously as in the case $\nu \neq 0$ for suitably chosen $D \subset Y \hookrightarrow Z$. 
	
	Now, we define the functionals $j,r \colon \BS \times \BS \to \setR$ and $b \colon I \times \BS \times \BS \to \setR$ by
	\begin{align}
		\bilJ(\chi, \phi)  &\coloneqq - \int_{\Omega} \div(\chi_2) \phi_1  \dx + \int_{\Omega} \chi_1 \div(\phi_2)\dx,  \\
		\bilR(\chi,\phi) & \coloneqq  \int_{\Omega} \gamma(\chi_1) \abs{\chi_2}^{s-2} \chi_2 \cdot \phi_2 \dx,\\
		\bilu(t,\chi,\phi) & \coloneqq - \int_{\partial \Omega } g(t)\, \phi_2 \cdot n \dsigma,\label{eq:forms-wave-b-v}
	\end{align} 	
	for any $\chi, \phi \in \BS$ and any $t \in I$. 
	To show that Assumption~\ref{as:problem-general}~\ref{itm:forms} is satisfied with $p = s$, we only have to argue that the estimate~\eqref{est:forms-wave-1} holds for $\BS$ as defined in~\eqref{def:VZ-wave-2}. 
	
	For $\bilJ$ the estimate still holds for any $\chi,\phi \in L^2(\Omega) \times H_{\Gamma}(\div;\Omega)$. 
	For $\bilR$ we have one term less than before, and the remaining term is estimated as before using the fact that the second component of functions in $\BS$ are contained in $L^s(\Omega)^d$. 
	Finally, for the bound on $\bilu(t,\chi,\phi)$ we use the fact that the trace operator is bounded from $H(\div;\Omega)$ to $H^{-1/2}(\partial \Omega)$, which is the dual space of the fractional space $H^{1/2}(\partial \Omega)$, see, e.g.,~\cite[Ch.~I.2.2]{Girault1986}. 
	Thus, it follows that 
	\begin{align*}
		\abs{b(t,\chi,\phi)} 
	&	\leq
		 \norm{g(t)}_{H^{1/2}(\partial \Omega)} \norm{\phi_2 \cdot n}_{H^{-1/2}(\partial \Omega)}\\
		& \leq c
		  \norm{g(t)}_{H^{1/2}(\partial \Omega)} \norm{\phi_2 }_{H(\div;\Omega)}\\
		  & \leq c
		  \norm{g(t)}_{H^{1}(\Omega)} \norm{\phi_2 }_{H(\div;\Omega)},
	\end{align*}
	and the remaining arguments are as above, assuming that $g \in C(I;H^{1}(\Omega))$. 
\end{enumerate}
\end{example}

\begin{example}[doubly nonlinear parabolic equation]\label{ex:doubly-parab}
	For functions $\alpha, \beta$ to be specified we consider the  nonlinear scalar evolution equation 
		\begin{alignat}{3} \label{eq:parab-doubly-eq}
			\partial_t \alpha(v) - \div (\beta(\nabla v)) &= f \qquad && \text{ on } Q, 
		\end{alignat}
	 for given function $f\colon Q \to \setR$. 
		This is supplemented by the following boundary conditions: 
	For some sufficiently smooth part of the boundary $\Gamma \subset \partial \Omega$, we impose 
	\begin{subequations}
		\begin{alignat}{3}\label{eq:parab-doubly-bc}
			v &= 0 \quad &&\text{ on } I \times \Gamma,\\
			\beta(\nabla v ) \cdot n +  \delta v & = g \qquad \qquad && \text{ on } I \times (\partial \Omega \setminus \Gamma),
		\end{alignat}
		\end{subequations}
 for a given function $g \colon I \times \partial \Omega \setminus \Gamma \to \setR$, and a constant $\delta>0$, with $n$ the outer unit normal on $\partial \Omega$. 
 Furthermore,  initial conditions are imposed.  

	The functions $\alpha\colon \setR \to \setR$ and $\beta \colon \setR^d \to \setR^d$ are assumed to be monotone and may be singular for $v = 0$, and for $\nabla v = 0 $, respectively. 
	In one space dimension $d = 1$  this equation is of  particular relevance to gas flow in pipelines in the high-friction and low Mach number regime, cf.~\cite{EggerGiesselmannKunkelEtAl2022}.	
	
	If $\alpha$ is invertible, then formally from~\eqref{eq:parab-doubly-eq} one may derive the equation for $z \coloneqq \alpha (v)$ as 
	\begin{align}\label{eq:parab-doubly-trafo}
		\partial_t z -  \div(\beta(\nabla(\alpha^{-1}(z)))) = f. 
	\end{align}
	In the following we assume that $\alpha$ and $\beta$ are given by
	\begin{subequations}\label{def:alpha-beta}
	\begin{alignat}{3} \label{def:alpha-pol}
		\alpha(v) &= \abs{v}^{-\sfrac{1}{q'}} v \qquad \;
		&&\text{ for some } q \in [1,\infty),\\ 
		\label{def:beta-pol}
		 \beta( \nabla v) &= \abs{\nabla v}^{p-2} \nabla v \qquad		&&\text{ for some } p \in (1,\infty), 
	\end{alignat}
	\end{subequations}
	with $q'\in (1,\infty]$ the Hölder conjugate to $q$, defined by $\frac{1}{q} + \frac{1}{q'} = 1$. 
	In the following we assume that $p\geq\frac{2d}{d+1}$.

For $\alpha, \beta$ as in~\eqref{def:alpha-beta} we obtain $v = \alpha^{-1}(z) = \abs{z}^{q-1} z$, and then~\eqref{eq:parab-doubly-trafo} reads
	\begin{align}\label{eq:parab-doubly-pq}
		\partial_t z -  \div\left(\abs{ \nabla v}^{p-2} \nabla v\right)
		= \partial_t z -  \div\left(\bigabs{ \nabla\left(\abs{z}^{q-1} z\right) }^{p-2} \nabla \left(\abs{z}^{q-1} z\right)\right) = f. 
	\end{align}
	Note that for $q = 1$ this evolution equation reduces to the \emph{$p$-Laplace equation}, and for $p = 2$ it reduces to the \emph{porous medium equation}. 
	
	Let us determine the setup for \eqref{eq:parab-doubly-pq} to fit into the framework described in~\eqref{eq:pH-weak} with Assumption~\ref{as:problem-general} satisfied. 
	For any $\phi \in C^{\infty}(Q)$, with $\phi|_{I \times \Gamma} = 0$, integrating by parts, employing the boundary conditions, and using $v = \alpha^{-1}(z)$ we obtain  
	\begin{align*}
		-	\int_{\Omega} \div \left( \beta(\nabla v)\right) \phi \dx 
		&= 
		-	\int_{\partial \Omega \setminus \Gamma}  \beta(\nabla v) \cdot n \, \phi \dx  + 
		\int_{\Omega}  \beta(\nabla v ) \cdot \nabla \phi \dx \\
		& = 
		\delta  \int_{\partial \Omega \setminus \Gamma}  v  \, \phi \dsigma  
		-  \int_{\partial \Omega \setminus \Gamma}  g \, \phi \dsigma 
		+ \int_{\Omega}  \beta(\nabla v ) \cdot \nabla \phi \dx \\
		& = 
		\delta  \int_{\partial \Omega \setminus \Gamma}  v  \, \phi \dsigma
		-  \int_{\partial \Omega \setminus \Gamma}  g \, \phi \dsigma 
		+ \int_{\Omega}  \abs{\nabla v}^{p-2} \nabla v  \cdot \nabla \phi \dx. 
	\end{align*} 
	Thanks to $p\geq \frac{2d}{d+1} >  \frac{2d}{d+2}$  we have a continuous and dense embedding 
	\begin{align}\label{def:VZ-parabol}
\BS \coloneqq W^{1,p}_{\Gamma}(\Omega) \hookrightarrow L^2(\Omega) \eqqcolon Z,
	\end{align}
	where $W^{1,p}_{\Gamma}(\Omega)$ denotes the closure of $\{f \in C^\infty(\Omega)\colon f|_{\Gamma} = 0\}$ with respect to the norm in $W^{1,p}(\Omega)$. 
	The Hamiltonian of~\eqref{eq:parab-doubly-pq} is given by 
	\begin{align*}
		\mathcal{H}(z) \coloneqq  \tfrac{1}{q+1} \int_{\Omega} \abs{z}^{q+1} \dx,
	\end{align*}
	which is defined for any $z\in  L^{q+1}(\Omega) \cap \BS$. 
	With the Fréchet derivative 	
	\begin{align*}
		\eta(z) \coloneqq \mathcal{H'}(z) = \alpha^{-1}(z) = \abs{z}^{q-1}z,
	\end{align*} 
	one can check, that $\mathcal{H}' \in C^1(D;\setR)$ for example for $D  \coloneqq Y \coloneqq L^{q+1}(\Omega)\cap X$.
	Thus, Assumption~\ref{as:problem-general}~\ref{itm:H} is satisfied, since $\eta(z) \in L^{(q+1)'}(\Omega) \subset Y'$ for any $z \in Y \subset L^{q+1}(\Omega)$. 

	With $\BS = W^{1,p}_{\Gamma}(\Omega)$ the functionals $\bilJ,\bilR \colon \BS\times \BS \rightarrow \mathbb{R} $ and $ \bilu \colon I \times \BS\times \BS \rightarrow \mathbb{R} $ are chosen as $\bilJ(v,\phi) \equiv 0$ and 
	\begin{align*}
		\bilR(v,\phi)& \coloneqq
		 \int_{\Omega} \abs{\nabla v}^{p-2} \nabla v \cdot \nabla \phi \dx + \delta \int_{\partial \Omega \setminus \Gamma} v \phi \dsigma, \\
		\bilu(t,v,\phi) &\coloneqq 
		\int_{\Omega} f(t) \,\phi \dx + \int_{\partial \Omega \setminus \Gamma} g(t)\, \phi \dsigma, 
	\end{align*}
	for any $v,\phi \in \BS$, and for given  $u \coloneqq (f,g) \in C(I;L^{2}(\Omega)) \times C(I;L^{2}(\partial \Omega))$.
	Obviously they are linear in $\phi$. 
	Due to $p \geq \frac{2d}{d+1}$, the trace operator is bounded from $W^{1,p}(\Omega)$ to $L^{2}(\partial \Omega)$.
	Hence, with a trace inequality, the embedding $W^{1,p}(\Omega)\hookrightarrow L^2(\Omega)$, as well as Hölder's and Young's inequalities one can show that 
	\begin{align}\label{est:forms-parab}
	\abs{r(v,\phi)} + \abs{b(t,v,\phi)} 
	&\leq 
	c  \left(1 + \norm{f(t)}_{L^2(\Omega)} +  \norm{g(t)}_{L^2(\partial \Omega)} +  \norm{v}_{W^{1,p}(\Omega)}^{\max(2,p)-1} 
	\right) \norm{\phi}_{W^{1,p}(\Omega)},
	\end{align}
	for any $v,\phi \in \BS$ and any $t\in I$. 
	Thus, the estimate in~\ref{itm:forms} is satisfied with exponent $\max(2,p)$. 

By linearity of $r$ and $b$ in $\phi$ and the estimate~\eqref{est:forms-parab} Lipschitz continuity in $\phi$ follows. 
Continuity in~$v$ is obvious for the linear term in $r$, and also holds for the nonlinear term of $r$ for $p\in (1,\infty)$. 
Since $b$ does not explicitly depend on $v$ the corresponding properties are trivial. 
To verify~\ref{itm:forms} we only have to show dissipativity of $r$.
Indeed, with $\delta \geq 0$ we have that
\begin{align*}
r(v,v)= \norm{\nabla v}_{L^p(\Omega)}^p  + \delta \norm{v}_{L^2(\partial \Omega)}^2 \geq 0 \quad \text{ for any } v \in \BS. 
\end{align*}
Note that the nonlinearity of $\bilR$ in $v$ stems from the $p$-Laplace structure, and the non-trivial $\eta(z)$ arises from the porous medium part of the equation.		

For well-posedness to~\eqref{eq:parab-doubly-eq} for $pq>1$ subject to homogeneous Dirichlet boundary conditions, i.e., $\Gamma = \partial \Omega$, see~\cite{Raviart1970}. For well-posedness in case $\Gamma = \emptyset$ see~\cite{SchoebelKroehn2020} for certain $p,q$.
\end{example}

\begin{remark}\label{rmk:power-extensions}
Several structural assumptions on $\alpha, \beta$ in Example~\ref{ex:doubly-parab} and on $F$ in Example~\ref{ex:wave-eq} can be relaxed.  
For example, alternative monotone functions, so-called Orlicz functions can be used, see, e.g.,~\cite{Diening2007}. 
Under the assumption of sufficiently regular functions, also vector-valued systems may be considered, but such regularity is rarely available. 
Also the regularity assumptions on the data can be somewhat weakened.
\end{remark}

\section{Structure-preserving time discretization}
\noindent 
In this section we introduce the general structure-preserving time-discrete numerical scheme for port-Hamiltonian systems of the form~\eqref{eq:pH-weak}. 
All of this equally applies  to the finite-dimensional case in~\eqref{eq:pH-findim}. 

In Section~\ref{sec:prelim} we collect the tools that are used in the following. 
Section~\ref{sec:scheme} introduces the time-discrete scheme with general quadrature. 
Furthermore, Proposition~\ref{prop:discr-energy} contains the conservation and dissipation of the Hamiltonian at the time grid points. 

\subsection{Preliminaries}\label{sec:prelim}
\noindent 
For a final time $T>0$ and a number $m \in \mathbb{N}$ we consider the collection of time points $\{t_0, \ldots, t_m\}$ with $0=t_0< t_1< \ldots < t_m = T$. 
 They generate a partition of $I\coloneqq[0,T]$ denoted by  
\begin{align*}
  I_\tau \coloneqq \{I_1, \ldots, I_m\},
\end{align*}
with subintervals $I_i \coloneqq [t_{i-1},t_i]$, for $i=1,\dots,m$, and denote the length of $I_i$ by $\tau_i \coloneqq t_{i} - t_{i-1}>0$. 
For convenience and with slight abuse of notation, we use the maximal mesh size $\tau \coloneqq\max_{i = 1, \ldots, m} \tau_i$ as index for the partition $I_\tau$. 

For a Banach space $\BS$ and for arbitrary $i=1,\dots,m$, by $\rP_k(I_i;\BS)$ we denote the set of polynomials of degree at most~$k \in \mathbb{N}_0$ mapping from $I_i$ to $\BS$. 
Furthermore, we define the semi-discrete (in time) function spaces of  piecewise polynomial functions of degree at most $k \in \mathbb{N}_0$ with values in $\BS$ as
\begin{align*}
	\rV_k(I_\tau;\BS) &\coloneqq \{ z \in L^\infty(I;\BS) \colon  z|_{I_i} \in \rP_k(I_i;\BS)\;  \text{ for all } i \in \{1, \ldots, m\} \}, \quad \text{and}\\ 
	\rV_k^c(I_\tau;\BS) &\coloneqq \rV_k(I_\tau;\BS)   \cap C(I; \BS).
\end{align*}

\subsubsection*{Quadrature.}
For the numerical scheme we require a quadrature rule to approximate the integrals on the intervals $I_i$, for $i \in \{1, \ldots, m\}$. 
For a general quadrature formula on $[0,1]$ with $s_{Q} \in \mathbb{N}$ nodes, the corresponding quadrature formulas on $I_i$ are $Q_i \colon C(I_i) \to \setR$ of the form
\begin{equation}\label{eq:quadrature-definition}
	Q_i(g)= \tau_i \sum_{j = 1}^{\quadnodes} \omega_j\, g(\zeta_j^i) \qquad \text{ for } i \in \{1, \ldots, m\},
\end{equation}
for some given weights $\omega_j $ and nodes $\zeta_j^i \in I_i$, $ i\in \{1, \ldots, m\}$, $j = 1, \ldots, \quadnodes$. 
For the weights we require that
\begin{align*}
	w_j >0 \quad \text{ for any } j \in \{1, \ldots, \quadnodes\} \qquad \text{ and} \qquad \sum_{j = 1}^{\quadnodes} w_j = 1.
\end{align*} 

\subsubsection*{$L^2$-projection.}
For a Hilbert space $Z$ the $L^2$-projection mapping to piecewise polynomial functions in time with values in $Z$ is a central tool in the numerical scheme proposed in the sequel. 
However, we also need to apply the projections to functions in $L^2(I;\BS)$ for a Banach space $\BS$.
Hence, as above we assume that $Z$ is a separable Hilbert space identified with its dual $Z'$, and that $\BS$ is a reflexive Banach space with continuous and dense embedding $\BS \hookrightarrow Z$, see Assumption~\ref{as:problem-general}.  
This means that $(\BS,Z,\BS')$ forms a Gelfand triple 
\begin{align*}
	\BS \hookrightarrow Z \cong Z' \hookrightarrow \BS',
\end{align*} 
with dense embedding $Z' \hookrightarrow \BS'$. 
With inner product $\skp{\cdot}{\cdot}$ on $Z$, we have in particular that 
\begin{align*}
	\skp{v}{w}_{\BS',\BS} = \skp{v}{w} \quad \text{ for any } v  \in Z, w \in \BS. 
\end{align*}

From now on let $k \in \mathbb{N}$ be fixed. 
The $L^2$-projection $\Pi \colon L^2(I;Z) \to \rV_{k-1}(I_\tau;Z)$ is defined by 
\begin{align}\label{def:L2proj}
\int_0^T \skp{(\Pi f) (t)}{g(t)} \dt = \int_0^T \skp{f(t)}{g(t)} \dt \qquad \text{ for all } g \in \rV_{k-1}(I_\tau;Z),
\end{align}
for $f \in L^2(I;Z)$. 
Note that since the functions in $\rV_{k-1}(I_\tau;Z)$ are discontinuous piecewise polynomials, the $L^2$-projection $\Pi$ is local. 
With the local $L^2$-projection $\Pi_i \colon L^2(I_i;Z)\to \rP_{k-1}(I_i;Z)$, this means that $\Pi f |_{I_i} = \Pi_i(f |_{I_i})$.  

For the sake of completeness let us summarize some stability properties of $\Pi$. 

\begin{lemma}[stability]\label{lem:Pi-stab}
	Let $k \in \mathbb{N}$ and let the function spaces as above.
Then the $L^2$-projection $\Pi \colon L^2(I;Z) \to \rV_{k-1}(I_\tau;Z)$ defined in~\eqref{def:L2proj} maps $L^2(I;\BS) \to \rV_{k-1}(I_\tau;\BS)$. 

 Furthermore, the following estimates are satisfied:
\begin{alignat*}{5}
\norm{\Pi f}_{L^p(I_i;Z)} &\leq c_p \norm{f}_{L^p(I_i;Z)}
\quad &&\text{ for any } f \in L^p(I_i;Z), \quad && \text{ for }  p \in [1,\infty], \\
\norm{\Pi f}_{L^p(I_i;\BS)} &\leq c_p \norm{f}_{L^p(I_i;\BS)}
\quad &&\text{ for any } f \in L^p(I_i;\BS),  \quad && \text{ for }  p \in (1,\infty].
\end{alignat*}
for any $i \in \{1, \ldots, m\}$. 
The constant $c_p>0$ depends only on $k$ and on $p$.  
The corresponding estimates hold for $I_i$ replaced by $I$. 
\end{lemma}
\begin{proof}We prove the statements on a single interval $I_i$ with $\Pi$ mapping to $\rP_{k-1}(I_i;Z)$, for $i \in \{1, \ldots, m\}$, since the global estimates are a consequence thereof. 
	
	The first estimate for $p = 2$ follows directly by definition of $\Pi$. 
	For general $p$ the estimate is proved by inverse estimates in time: 
	Indeed, for $p\geq 2$ applying an inverse estimate in combination with the stability in $L^2(I_i;Z)$ as well as Hölder's inequality we obtain
	\begin{align*}
		\norm{\Pi f }_{L^p(I_i;Z)} 
		\leq c \tau_i^{\frac{1}{p} - \frac{1}{2}}
		\norm{\Pi f }_{L^2(I_i;Z)} 
		\leq c \tau_i^{\frac{1}{p} - \frac{1}{2}}
		\norm{ f }_{L^2(I_i;Z)} 
		\leq c 
		\norm{ f }_{L^p(I_i;Z)}. 
	\end{align*} 
	The constant depends only on $k \in \mathbb{N}$ and on $p$. 
	
	To prove the statement for $p<2$, we use the fact that $(L^{p}(I_i;Z))' = L^{p'}(I_i;Z)$.
	We obtain by duality that
	\begin{align*}
		\norm{\Pi f}_{L^p(I_i;Z)}  
		&=
		 \sup_{g \in L^{p'}(I_i;Z)} \frac{\int_{I_i} \skp{\Pi f(t)}{g(t)} \dt }{\norm{g}_{L^{p'}(I_i;Z)}} 
	 =
		 \sup_{g \in L^{p'}(I_i;Z)}  \frac{\int_{I_i} \skp{ f(t)}{\Pi g(t)} \dt }{\norm{g}_{L^{p'}(I_i;Z)}} 
		  \\ 
		  &\leq \norm{f}_{L^p(I_i;Z)} \sup_{g \in L^{p'}(I_i;Z)} \frac{\norm{\Pi g}_{L^{p'}(I_i;Z)}}{\norm{g}_{L^{p'}(I_i;Z)}}
		  \leq c \norm{f}_{L^p(I_i;Z)}. 
	\end{align*} 
In the last step we have used the previously proved stability of $\Pi$ in $L^{p'}(I_i;Z)$ with $p'\geq 2$. 

	To prove the second identity let $f \in L^p(I_i;\BS) \hookrightarrow L^p(I_i;Z)$ for $p \in (1,\infty)$ be arbitrary. 
	Then the linear functional represented by $\Pi f$ satisfies 
	\begin{align*}
		\ell(g) \coloneqq \bigabs{ \int_{I_i} \skp{\Pi f}{g} \dt }  = \bigabs{ \int_{I_i} \skp{f}{g} \dt }\leq \norm{f}_{L^p(I_i;\BS)} \norm{g}_{L^{p'}(I_i;\BS')} 
	\end{align*}
	for any $g \in L^2(I_i;Z)$. 
	Consequently, with the (dense) embedding $L^{p'}(I_i;Z) \hookrightarrow L^{p'}(I_i;\BS')$ the Hahn--Banach extension theorem shows that the bounded linear functional $\ell \colon \rP_{k-1}(I_i;Z)\to \setR$ uniquely extends  to a bounded linear functional  $\ell \colon \rP_{k-1}(I_i;\BS') \to \setR$ with the same operator norm. 
	By the uniqueness of the extension and the fact that $L^p(I_i;\BS) \subset L^q(I_i;\BS)$ for $p \geq q$ the extension is  unique and independent of $p$. 
	In particular, since $\BS$ is reflexive, we have that $\Pi f \in \rP_{k-1}(I_i;\BS)$ and that
	\begin{align*}
		\norm{\Pi f}_{L^p(I_i;\BS)}  \leq c \norm{f}_{L^p(I_i;\BS)},
	\end{align*}
	with constant independent of $I_i$ and $f$. 
	The estimate for $p = \infty$ follows again by an inverse estimate. 
	This proves the second estimate.  
\end{proof}

\subsection{Petrov--Galerkin scheme}\label{sec:scheme}
\noindent
In this section we present the structure-preserving scheme to approximate solutions to the port-Hamiltonian system~\eqref{eq:pH-weak} in the setting of Assumption~\ref{as:problem-general}. 
It is a Petrov--Galerkin type approximation, meaning that it uses a variational formulation with different polynomial degrees for the solution and for the space of test functions. 
More specifically, it is a \emph{continuous Petrov--Galerkin method}, cf.~\cite[Sec.~70.1.2]{Ern2021}. 
Here 'continuous' refers to the fact that the trial space consists of continuous piecewise polynomial functions in time. 
We introduce it in combination with some quadrature formula $Q_i$ approximating the integral on $I_i$ for all nonlinear terms, see \eqref{eq:quadrature-definition}.

\begin{scheme}\label{def:scheme}
	Find $\zh \in \rV_k^c(I_\tau; \BS)$ such that $\zh(0) = z_0$ and
	\begin{equation}\label{eq:scheme_global}
		\int_0^T \skp{ \del_t \zh }{ \phi} \dt
		= \sum_{i = 1}^m 
		Q_i \left[ \bilJ(\Pi \eta(\zh), \phi) 
		- \bilR(\Pi \eta(\zh), \phi) 
		+ \bilu(\cdot,\Pi\eta(\zh),\phi) \right] 
	\end{equation}
	holds for all $\phi \in \rV_{k-1}(I_\tau; \BS)$.
	
	Since the test functions are discontinuous,~\eqref{eq:scheme_global} can be localized in time. Indeed, we can equivalently reformulate it as a time stepping via 
	\begin{equation}\label{eq:scheme_local}
		\int_{t_{i-1}}^{t_i} \skp{ \del_t \zh }{ \phi} \dt
		= 
	 Q_i\left[\bilJ(\Pi \eta(\zh), \phi) 
		- \bilR(\Pi \eta(\zh), \phi) 
		+ \bilu (\cdot,\Pi\eta(\zh),\phi) \right],  
	\end{equation}
	for all $\phi \in \rP_{k-1}(I_i;\BS)$ and all $i=1,\dots,m$.
\end{scheme}

\newpage

\begin{remark}\label{rmk:scheme}\hfill
\begin{enumerate}[label = (\alph*)]
	\item The $L^2$-projection in the term $b$ is not essential, and it has no relevance in case $b$ is state-independent. 
	\item In the special case that $j,r$ are bilinear and that $\eta$ is the identity, i.e., $\mathcal{H}$ is quadratic, the $L^2$-projection in~\eqref{eq:scheme_global} vanishes. 
Therefore, it is straightforward to see, that the method reduces to the standard continuous Petrov--Galerkin method in this case. 
In particular, for $k = 1$ and midpoint quadrature rule the scheme reduces  to the implicit midpoint method.  
\item 
Let us briefly discuss under which conditions the terms in~\eqref{eq:scheme_global} are well-defined. 	
For this purpose we assume that a discrete solution $\zh \in \rV_{k}^c(I_\tau;\BS) \subset  C(I;\BS)$ satisfies that $\eta(\zh) \in C(I;\BS)$.  
Consequently, by Lemma~\ref{lem:Pi-stab} we have that $\Pi \eta(\zh) \in \rV_{k-1} (I_\tau;\BS)$, i.e., in particular it is piecewise continuous. 
The same is true for $\phi$, and with continuity of $\bilJ$ in both arguments, it follows that $\bilJ(\Pi \eta(\zh),\phi)$ is piecewise continuous in time with respect to $I_\tau$.  
Therefore, the quadrature is well-defined, and the arguments of the other terms on the right-hand side proceed analogously. 
The left-hand side is well-defined thanks to $\partial_t \zh, \phi \in \rV_{k-1}(I_\tau;\BS)$ and the fact that $\BS \hookrightarrow Z$. 

The scheme can also be considered without the quadrature, by integrating over $I_i$ rather than using a quadrature $Q_i$. 
Integrability of all terms follows by the above arguments, by the estimates in Assumption~\ref{as:problem-general}~\ref{itm:forms}, and by the stability properties of~$\Pi$ in Lemma~\ref{lem:Pi-stab}. 
\item  
For practical purposes the $L^2$-projection has to be computed using a further quadrature formula to approximate $\int_{I_i} \skp{\eta(\zh)}{g} \dt$ for any $\phi \in \mathcal{P}_{k-1}(I_i;Z)$, cf.~\eqref{def:L2proj}. 
\end{enumerate}
\end{remark}

\begin{remark}[comparison of schemes]\label{rmk:discr-comp}
	The scheme in \cite[Sec.~7]{morandin24-modeling} for infinite-dimensional and finite-dimensional systems of the form \eqref{eq:pH-findim-classical} reads: 
	Find $\zh \in \rV_k^c(I_\tau; \BS)$ such that $\zh(0) = z_0$ and 
\begin{equation}\label{eq:scheme_global-v}
	\int_0^T \skp{ \del_t \zh }{ \phi} \dt
	= \int_{0}^T 
	\left[ (\tilde J(\zh) \Pi \eta(\zh), \phi) 
	- (\tilde R(\zh) \Pi \eta(\zh) , \phi) 
	+ (\tilde B(\zh)u, \phi)  \right] 
\end{equation}
holds for all $\phi \in \rV_{k-1}(I_\tau; \BS)$. 
For singular $\tilde J$ and $\tilde R$ difficulties may occur since at some points in space one may have that, e.g., $\zh = 0$ and $\Pi(\eta(\zh)) \neq 0$. 
Singular operators include the friction term~\eqref{eq:friction-powerlaw} for $s<2$ in the quasilinear wave equation in Example~\ref{ex:wave-eq} and operator $R$ for the $p$-Laplace equation for $p < 2$, see Example~\ref{ex:doubly-parab}.  
For the latter example there is a singularity in the gradient. 
 In~\eqref{eq:scheme_global-v} the term $\tilde R(\zh) \Pi \eta(\zh) = \abs{\nabla \zh}^{p-2} \Pi (\nabla \zh)$ would occur. 
In contrast, the corresponding term in Scheme~\ref{def:scheme} reads
\begin{align*} 
r(\Pi \eta(\zh),\phi) = \int_{\Omega} \abs{\nabla \Pi \zh}^{p-2} \nabla \Pi \zh \cdot \nabla \phi \dx,
	\end{align*}
which is monotone and well-defined also for $\nabla \zh = 0$. 
\end{remark}

By construction, sufficiently smooth solutions to the scheme satisfy a discrete version of the energy balance~\eqref{eq:energybalance}, which we refer to as \emph{energy consistency}. 

\begin{proposition}\label{prop:discr-energy}
	Solutions 
	$\zh \in \rV_{k}^c(I_\tau;\BS)$ 
	 to Scheme~\ref{def:scheme} 
	 with $\eta(\zh) \in C(I;X)$
	  satisfy 
	\begin{equation}\label{eq:energybalance-projecte}
		\begin{aligned}
		\cH(\zh(t_i)) - \cH(\zh(t_{i-1})) &= Q_i\left[
		-  \bilR(\Pi\eta(\zh),\Pi \eta(\zh)) 
		+ \bilu( \cdot,\Pi \eta(\zh), \Pi \eta(\zh)) 
		\right] \\
		& \leq  Q_i \left[ 
		\bilu (\cdot,\Pi \eta(\zh),\Pi \eta(\zh)) \right],
		\end{aligned}
	\end{equation}
	for all $i = 1,\dots,m$.
\end{proposition}
\begin{proof}
 Let $i \in \{1,\dots, m\}$ be arbitrary and let $\zh \in \rV_{k}^c(I_\tau;\BS)$ with $\eta(\zh) \in  C(I;\BS)$ be a solution to Scheme~\ref{def:scheme}. 
 Then, we have in particular $\partial_t \zh \in \rV_{k-1}(I_\tau;\BS) \subset L^\infty(I;\BS)$. 
Since $\Pi$ is the $L^2$-orthogonal projection to $\rV_{k-1}(I_\tau; Z)$  using~\eqref{eq:scheme_local} we find that
\begin{align*}
	\cH(\zh(t_i)) &- \cH(\zh(t_{i-1}))
	= \int_{t_{i-1}}^{t_i}\frac{\mathrm{d}}{\mathrm{d}t} \cH(\zh) \dt 
	= \int_{t_{i-1}}^{t_i} \langle \cH'(\zh), \del_t \zh\rangle \dt \\
	&= \int_{t_{i-1}}^{t_i} \langle \eta(\zh),  \del_t \zh\rangle \dt
	 = \int_{t_{i-1}}^{t_i}\langle \Pi \eta(\zh), \del_t \zh\rangle \dt \\
	& = Q_i \left[ \bilJ(\Pi \eta(\zh), \Pi \eta(\zh)) 
	- \bilR(\Pi \eta(\zh), \Pi \eta(\zh)) 
	+ \bilu (\cdot, \Pi \eta(\zh),\Pi \eta(\zh)) \right]. 
	\end{align*}
	Employing the conservation property of $\bilJ$ and the dissipative nature of $\bilR$ due to Assumption~\ref{as:problem-general}~\ref{itm:forms}, pointwise at the quadrature nodes, yields  
	\begin{align*}
	\cH(\zh(t_i)) - \cH(\zh(t_{i-1}))
	& = Q_i \left[  
	- \bilR(\Pi \eta(\zh), \Pi \eta(\zh)) 
	+ \bilu (\cdot,\Pi \eta(\zh),\Pi \eta(\zh))\right] \\
	& \leq Q_i \left[ \bilu (\cdot,\eta(\zh),\Pi \eta(\zh))\right], 
\end{align*}
which proves the claim. 
\end{proof}

\begin{remark}\label{rmk:comp-energyconsist}\hfill
	\begin{enumerate}[label = (\alph*)]
		\item 
		For finite-dimensional port-Hamiltonian systems of the form~\eqref{eq:pH-findim} the discrete energy balance reduces to 
		\begin{align*} 
				\cH(\zh(t_i)) - \cH(\zh(t_{i-1}))
				& = Q_i \left[  
				- \skp{R(\Pi \eta(\zh))}{ \Pi \eta(\zh)}_{\ell^2}
				+  \skp{B(\cdot,\Pi\eta(\zh))}{\Pi \eta(\zh)}_{\ell^2}\right]\\
				& \leq Q_i \left[  
			 \skp{B(\cdot,\Pi\eta(\zh))}{\Pi \eta(\zh)}_{\ell^2} \right],
		\end{align*}
		and the discrete output can be defined from this. 
		\item Note that for the special case of a Hamiltonian system, i.e., with $r \equiv b \equiv 0$, the Hamiltonian is exactly conserved by our scheme. 
		This is a stronger property than, e.g., the one achieved in~\cite{mehrmann19-structure} for general Hamiltonian functions. 
		There, the estimate is satisfied only asymptotically for $\tau \to 0$. 
		\item 	Proposition~\ref{prop:discr-energy} makes a statement for $\mathcal{H}(\zh(t))$ for $t = t_i$, but not for arbitrary time points. 
		In case $\bilR$ is coercive in all components,  boundedness of $\mathcal{H}(\zh(t))$ can be obtained for any $t \in (t_{i-1}, t_i)$ using methods presented in~\cite{ChrysafinosWalkington2010}. 
		Since those arguments are available only for coercive $\bilR$ we refrain from presenting the details. 
		\item For the scheme in~ \cite{morandin24-modeling} a corresponding discrete energy balance is available. 
	\end{enumerate}
\end{remark}


\section{Numerical Experiments}\label{sec:num-exp}

\noindent To showcase the performance of our cPG scheme presented in Section~\ref{sec:scheme} we test it on some of the examples in Section~\ref{sec:ex}. 
This includes finite-dimensional port-Hamiltonian systems as well as space-discrete versions of infinite-dimensional ones. 
For polynomial degree $k \in \mathbb{N}$ in Scheme~\ref{def:scheme} we use Gauß quadrature with $s_Q \in \mathbb{N}$ nodes for $Q_i$ and Gauß quadrature with $s_\Pi \in \mathbb{N}$ nodes to approximate $\Pi(\eta(z_\tau))$, cf.~\eqref{eq:scheme_local}. 
We numerically investigate the impact of the choice of $s_Q$ and $s_\Pi$ on convergence and energy consistency. 
  
We observe  optimal rates $\tau^{k+1}$ for sufficiently smooth solutions, as for the standard cPG method, cf.~\cite{Aziz1989}, provided that $\quadnodes, \projnodes \geq k$. 
This means that the use of the $L^2$-projection does not affect the convergence rate and neither does its approximation by means of quadrature. 
Furthermore, under the same conditions we observe superconvergence of order $\tau^{2k}$ at the time grid points, as expected from~\cite{Aziz1989}. 

To investigate the energy consistency recall that the property in Prop.~\ref{prop:discr-energy} is proved with quadrature $Q_i$, but without use of quadrature for computing $\Pi \eta(\zh)$. 
Consequently, this property is independent of the choice of $s_Q$, but the effect of the choice of $s_\Pi$ on the energy consistency has to be addressed. 
In the examples under consideration we find that using quadrature on the $L^2$-projection energy consistency is not satisfied exactly unless $\eta$ is polynomial, and $\projnodes$ is chosen such that exactness is ensured. 
However, for $s_\Pi$ sufficiently large the relative error is close to machine precision. 

In Section~\ref{sec:space-discr} we summarize a structure preserving space discretization, similar to the one in \cite{Jackaman2019, JackamanPryer2021}, and \cite[Sec.~7.2]{morandin24-modeling}.

\subsection{Implementation details}
\noindent 
For finite-dimensional  and space-discrete port-Hamiltonian systems we consider examples of \eqref{eq:pH-findim} with $X = Z=\mathbb{R}^{\delta}$, where $\delta \in \mathbb{N}$ is the dimension of the system, i.e., in the latter case the number of degrees of freedom of the space discretization. 
To such systems we apply Scheme~\ref{def:scheme}. 
Since it localizes as described above, in each time step one has to solve for $\zh|_{I_i}$ with imposed value for $\zh(t_{i-1})$. 
For example, $\zh$ can be expanded in the $L^2$-normalized Legendre polynomials forming a basis of $\mathcal{P}_{k}(I_i)$. 
Note that for $s_\Pi = k$ quadrature nodes and if $\eta = \identity$, then the integrals in the $L^2$-projection are evaluated exactly. 

In each time step  a nonlinear system of equations has to be solved, for which we use Newton's method. 
The derivatives required in Newton's method are computed symbolically using \texttt{JAX}~\cite{bradbury18-jax}. 
In the first time step the constant one vector is used as starting value for the Newton iteration, and in subsequent time steps the numerical solution of the previous time step is used. 
We have not encountered any issues regarding convergence of the Newton iteration in any of the computations. 

To study the convergence of the scheme we employ manufactured solutions. 
This means that for a function $z$ and a system determined by $J,R,\eta$ and $B$ we compute $\overline{B}$ and $\overline z_0$ such that 
\begin{equation}\label{eq:manuf-sol}
			\partial_t z
			= J(\eta(z))- R(
			\eta(z))
			+ \overline{B}(\cdot,\eta(z))
	\end{equation}
subject to $z(0) = \overline z_0$
is satisfied. 
This system is again of the form~\eqref{eq:pH-findim} and its exact solution $z$ is available. 

In the following we approximate the $L^\infty$-norm of the error $z - z_\tau$ by evaluation on  a time grid $\mathcal{T}_{\mathrm{fine}}  \coloneqq \{ i \tauref\colon i \in \{0, ..., 38913\}\}$ with time step size $\tauref \approx 1.29 \cdot 10^{-4}$. The corresponding relative error is denoted by 
	\begin{align}\label{def:error-E}
		E \coloneqq \frac{\max_{t \in \mathcal{T}_{\mathrm{fine}}} \norm{z_{\tau}(t) - z(t)}}{\max_{t \in \mathcal{T}_{\mathrm{fine}}} \norm{ z(t)}}.
	\end{align}
Furthermore, we evaluate the error at the time grid points to investigate nodal superconvergence
\begin{align}\label{def:error-Etau}
	E_{\tau} \coloneqq \frac{\max_{t \in \{t_0,t_1 , \ldots, t_m\}} \norm{z_{\tau} (t)- z(t)}}{\max_{t \in \{t_0,t_1, \ldots, t_m\}} \norm{ z(t)}}. 
\end{align}
To verify the proposed energy consistency of the method we compute the quantity 
\begin{equation}\label{def:E}
	\mathcal{E}(z_\tau; t_i) \coloneqq \frac{\abs{ \mathcal{H}(z_\tau(t_{i})) - \mathcal{H}(z_\tau(t_{i-1})) - Q_i[ - r(\widetilde \Pi\eta(z_\tau),\widetilde \Pi\eta(z_\tau)) + b(\cdot,\widetilde \Pi\eta(z_\tau), \widetilde \Pi\eta(z_\tau))] } }{ \max_{j=1,\dots,m} \abs{ \mathcal{H}(z_\tau(t_{j})) - \mathcal{H}(z_\tau(t_{j-1})) } },
\end{equation}
where $\widetilde \Pi$ denotes the approximation of the $L^2$-projection obtained by using quadrature. 
This represents a measure for the error in the energy balance in~\eqref{eq:energybalance-projecte}.

\subsection{Toda lattice}
\noindent 
To approximate solutions to Example~\ref{ex:toda} we use $N=5$, $\gamma_i = 0.1$ for $i =1,\dots,N$ and we numerically approximate the discrete solutions on the time interval $[0,T] = [0,5]$ with control input $u(t) = \sin(2t)$. 

To investigate convergence we choose the manufactured solution 
\begin{equation*}
  q_i(t) = \sin(t) \quad \text{ and } \quad  p_i(t) = \cos(t), \quad i = 1,\dots,N
\end{equation*}
and compute the corresponding term $\overline{B}$ and $\overline z_0$ in~\eqref{eq:manuf-sol}. 
In order to approximate the error between the approximate solution $\zh$ and the exact solution $z$ in $L^\infty(I)$ we  evaluated the difference $z - \zh$ on a time grid with step size~$\tauref \approx 1.29 \cdot 10^{-4}$. 
Figure~\ref{fig:toda_varying_degree} shows the convergence for polynomial degrees~$k \in \{1,2,3,4\}$ using Gauß quadrature with $s_Q = s_\Pi = k$ quadrature nodes.
Evidently, our method achieves the optimal convergence rate $\tau^{k+1}$. 
Figure~\ref{fig:toda_varying_degree_different_sampling} uses the same settings, but computes the errors only at the time grid points $t_0,\dots,t_m$.
Here, we observe convergence rate $\tau^{2k}$, i.e., our method exhibits nodal superconvergence. 

To verify that $s_Q = s_\Pi = k$ is the best choice for convergence  Figure~\ref{fig:toda_varying_quadrature} shows the 
convergence plot for polynomial degree $k=3 = s_\Pi$ and Gauß quadrature rules $Q_i$ with varying $s_Q$.  
Indeed, as expected the results illustrate that higher order quadrature rules do not improve the convergence rates and lower order quadrature rules reduce it.
Similarly, we investigate the effect of the choice of the quadrature used to compute the $L^2$-projection. 
For polynomial degree $k = 3 = s_Q$,  different Gauß quadrature rules with $s_\Pi$ nodes are employed for the approximation of the projection $\proj\eta(\zh)$. 
The results in  Figure~\ref{fig:toda_varying_projection} demonstrate that using quadrature rules with order of exactness higher than $2k-1$ does not improve the convergence, whereas quadrature rules with lower exactness degree lead to a reduced convergence rate.
Since in this example $\eta$ is nonlinear, this is a meaningful case to test the impact of the quadrature rule used in the computation of the $L^2$-projection.

Finally, we visualize the energy consistency property of our scheme by considering~\eqref{eq:toda} for $u(t) = \sin(2t)$ and $z_0 = 0$. 
In Figure~\ref{fig:toda_energybalance} the quantity $\mathcal{E}(\zh;t_i)$ in \eqref{def:E}, which is related to the relative error in the energy balance~\eqref{eq:energybalance-projecte}, is plotted for $\tau = 10^{-2}$ and  several polynomial degrees~$k\in\{1,2,3,4\}$, $s_Q = k$ and several numbers of Gauß quadrature nodes $\projnodes$ in the projection~$\Pi \eta(\zh)$ over~$t_i$. 
The experiment shows that our method satisfies the energy balance close to machine precision, provided that $\projnodes$ is sufficiently large. 
More specifically, one can observe that $s_\Pi \geq \max(k,3)$ yields satisfactory energy consistency. 

\begin{figure}
    \hfill
	\begin{subfigure}[t]{.3\textwidth}
		\centering
		\scalebox{0.31}{\input{figures/toda_varying_degree.pgf}}
		\caption{
			convergence in $\tau$ for 
		 several $k = \quadnodes = \projnodes $
		}
		\label{fig:toda_varying_degree}
	\end{subfigure}
	\hfill
	\begin{subfigure}[t]{.3\textwidth}
		\centering
		\scalebox{0.31}{\input{figures/toda_varying_quadrature.pgf}}
		\caption{
			convergence in $\tau$ for 
			 $k  = \projnodes = 3 $ and several $s_Q$}
		\label{fig:toda_varying_quadrature}
	\end{subfigure}
	\hfill 
	\begin{subfigure}[t]{.3\textwidth}
		\centering
		\scalebox{0.31}{\input{figures/toda_varying_projection.pgf}}
		\caption{
			convergence in $\tau$ for 
			$k  = \quadnodes = 3 $ and several $\projnodes$}
		\label{fig:toda_varying_projection}
	\end{subfigure}
	\hfill
	\caption{
		Convergence in $\tau$ for the Toda lattice~\eqref{eq:toda} for several polynomial degrees $k$ and several values of $\quadnodes$ and $\projnodes$.}
\end{figure}

\begin{figure}
    \centering
   	\begin{subfigure}[t]{.3\textwidth}
    	\centering
    	\scalebox{0.31}{\input{figures/toda_varying_degree_different_sampling.pgf}}
    	\caption{
    	nodal superconvergence for 
    	several $k = \quadnodes = \projnodes$
    	}
    	\label{fig:toda_varying_degree_different_sampling}
    \end{subfigure}
	\hspace{0.7em}
	\begin{subfigure}[t]{.3\textwidth}
		\centering
		\scalebox{0.31}{\input{figures/toda_energybalance.pgf}}
		\caption{
			relative energy error~\eqref{def:E}, $\tau = 10^{-2}$,  for
		 several $k = \quadnodes $ and $\projnodes$}
		\label{fig:toda_energybalance} 
	\end{subfigure}
	\caption{
		Nodal superconvergence in $\tau$ and relative error in the energy balance  for the Toda lattice~\eqref{eq:toda} for several polynomial degrees $k$ and several values of $\quadnodes$  and $\projnodes$. 
	}
\end{figure}

\subsection{Spinning rigid body}
\noindent Let us consider Example~\ref{ex:spinning-body} for given $Q$ and $\tilde{B}$ with $I_i = b_i = 1$, for $i =1,\dots,3$ and for control input  $u(t) = \sin(2t)$ on the time interval $[0,T] = [0,5]$. 
 Note that there is no dissipation since $R \equiv 0$. 
We compute the numerical solution $\zh$ with Scheme~\ref{def:scheme} for several time steps $\tau$ using Gauß quadrature with $s_Q = k$ quadrature nodes for $Q_i$ and $s_\Pi = k$ nodes for the approximation of $\Pi \eta(\zh)$. 
Since the Hamiltonian $\mathcal{H}$ is quadratic in this example, and hence $\eta$ is linear, this means that the computation of the $L^2$-projection is exact when using $s_\Pi = k$. 
Hence, by Proposition~\ref{prop:discr-energy} the energy at the time grid points is exactly preserved. 
  
As manufactured solution we choose $s = (p_1,p_2,p_3)\tp$ with
\begin{equation*}
  p_1(t) = \sin(t), \qquad p_2(t) = \sin(2t)\cos(t)^2 + 0.5, \qquad \text{ and }\quad  p_3(t)=\cos(t). 
\end{equation*}
To investigate the convergence we approximate the $L^\infty(I)$-norm of the error by $E$ as in \eqref{def:error-E} for~$\tauref \approx 1.29 \cdot 10^{-4}$ and $E_\tau$ as defined in \eqref{def:error-Etau}. 
Figure~\ref{fig:rigid_body_varying_degree} shows the convergence for polynomial degrees~$k \in \{1,2,3,4\}$ and $\quadnodes = \projnodes = k$. 
Again, we observe optimal convergence rate $\tau^{k+1}$.
Figure~\ref{fig:rigid_body_varying_degree_different_sampling} uses the same parameters, but displays the errors $E_\tau$ at the time grid points $t_0,\dots,t_m$  only. 
As before, our method exhibits nodal superconvergence with convergence rate $\tau^{2k}$.

Again, we visualize the energy consistency of our scheme by considering~\eqref{eq:rigid-body} for $Q$, $\tilde B$ and $u$ as above and $z_0 = (0, 0.5, 1)\tp$. 
In Figure~\ref{fig:rigid_body_energybalance} the  quantity $\mathcal{E}(\zh;t_i)$ as defined in~\eqref{def:E} is plotted over $t_i$, for $\tau = 10^{-2}$ and several polynomial degrees $k\in\{1,2,3,4\}$, $s_Q = k$ and $\projnodes = k$ Gauß quadrature nodes in the projection $\Pi \eta(\zh)$. 
The experiment shows that the energy balance is satisfied up to machine precision as expected. 

Note that since $\tilde J$ is linear in $z$ it is possible to choose the quadrature $Q_i$ in a way that the integration is exact. 
For this purpose it suffices to choose $\quadnodes \geq \frac{3}{2}k-1$.

\begin{figure}
    \hfill
	\begin{subfigure}[t]{.3\textwidth}
		\centering
		\scalebox{0.31}{\input{figures/rigid_body_varying_degree.pgf}}
		\caption{
			convergence in~$\tau$ for several~$k = \quadnodes = \projnodes$}
		\label{fig:rigid_body_varying_degree}
	\end{subfigure}
	\hfill
	\begin{subfigure}[t]{.3\textwidth}
		\centering
		\scalebox{0.31}{\input{figures/rigid_body_varying_degree_different_sampling.pgf}}
		\caption{
			nodal superconvergence in~$\tau$ for several~$k = \quadnodes = \projnodes$}
		\label{fig:rigid_body_varying_degree_different_sampling}
	\end{subfigure}
	\hfill
	\begin{subfigure}[t]{.3\textwidth}
		\centering
		\scalebox{0.31}{\input{figures/rigid_body_energybalance.pgf}}
		\caption{
			relative energy error \eqref{def:E}, $\tau = 10^{-2}$, for several $k = \quadnodes =\projnodes$}
		\label{fig:rigid_body_energybalance} 
	\end{subfigure}
		\hfill 
	\caption{
		Convergence~\eqref{fig:rigid_body_varying_degree}, 
		nodal superconvergence \eqref{fig:rigid_body_varying_degree_different_sampling} in $\tau$ and 
		relative error in the energy balance~\eqref{fig:rigid_body_energybalance} for the rigid spinning body~\eqref{eq:rigid-body}
		for several polynomial degrees $k$ and $\quadnodes = \projnodes = k$. 	 	
	}
\end{figure}

\subsection{Space discretization}\label{sec:space-discr}
\noindent 
Even though the discretization in space is not the primary focus of this work let us present a simple approach that yields a discrete Hamiltonian structure and that is closely linked to our approach for the time discretization, see also~\cite{Jackaman2019,JackamanPryer2021}. 
We start from a port-Hamiltonian system of the form~\eqref{eq:pH-weak} for $z\colon I \to \BS $ with $\BS$ a Banach space on a domain $\Omega$. 
We assume that $\BS_h\subset \BS$ is a finite dimensional subspace with a fixed basis with suitable approximation properties. 
This means that each $z_h \in \BS_h$ is determined by its coefficient vector $w_h \in \setR^{\dim \BS_h}$ in this basis and by $I_h$ we denote the canonical interpolation operator with $I_h(w_h) = z_h$. 
Furthermore, by $c_h$ we denote the mapping $z_h \mapsto c_h(z_h)  \coloneqq w_h$ of a function in $\BS_h$ to its coefficients in the chosen basis. 

\subsubsection*{Space-discrete port-Hamiltonian structure} 
To introduce the (space)-discrete Hamiltonian structure an orthogonal projection is instrumental again. 
Recalling that $\BS_h \subset \BS \hookrightarrow Z$ for a Hilbert space $Z$, let $\Pi_h \colon Z \to \BS_h$ be the $Z$-orthogonal projection.  
Now, let us define the discrete Hamiltonian $\mathcal{H}_h \colon \mathbb{R}^{\dim \BS_h}\to \mathbb{R}$ by
\begin{align*}
	\mathcal{H}_h(w_h) \coloneqq \mathcal{H}(I_h w_h). 
\end{align*}
Then, by linearity of $I_h$ and the projection property of $\Pi_h$ the corresponding $\eta_h$ with $\eta_h = \mathcal{H}_h'$ can be computed as follows for $v_h, w_h \in \setR^{\dim \BS_h}$
\begin{align*}
	\frac{\mathrm{d}}{\mathrm{d}s}\mathcal{H}_h(w_h + s v_h)\big|_{s = 0} 
	&= 
		\frac{\mathrm{d}}{\mathrm{d}s}\mathcal{H}(I_h (w_h) + s I_h(v_h))|_{s = 0} = \skp{\mathcal{H}'(I_h(w_h))}{I_h(v_h)} \\ 
		&= \skp{\eta(I_h(w_h))}{I_h(v_h)} 
		= \skp{\Pi_h \eta(I_h(w_h))}{I_h(v_h)} \\
		& = \skp{M c_h(\Pi_h \eta(I_h(w_h)))}{ v_h}_{\ell^2}
		\eqqcolon \skp{\eta_h(w_h)}{v_h}_{\ell^2}. 
\end{align*}
Here $M$ is the corresponding mass matrix of the basis of $\BS_h$ with respect to the inner product of $Z$. 
Note that this fully determines $\eta_h \colon \setR^{\dim \BS_h} \to \setR^{\dim \BS_h}$.  
This motivates the following space-discrete scheme
\begin{subequations}\label{eq:pH-space-discrete}
\begin{align}
\skp{\partial_t z_h}{\phi_h}  
&=  j(\Pi_h \eta(z_h),\phi_h) - r(\Pi_h \eta(z_h),\phi_h) + b(\cdot,\Pi_h \eta(z_h),\phi_h) \qquad \text{ for all } \phi_h \in \BS_h,
\\
z_h(0) &= \Pi_h(z_0). 
 \end{align}
 \end{subequations}
 This still has the form~\eqref{eq:pH-weak} replacing $j(\eta(z_h),\phi_h)$ by $j(\Pi_h \eta(z_h),\phi_h)= j(\Pi_h \eta(z_h),\Pi_h \phi_h)$, and hence the property~\ref{itm:bJ} is preserved. 
 For this reason, by the arguments above or by direct calculation energy consistency is available for any space-discrete solution $z_h$ as
 \begin{equation}\label{est:energy-consistency-h}
 \begin{aligned}
\mathcal{H}(z_h)(t) - \mathcal{H}(z_h)(t_0) &=  \int_{t_0}^t - r(\Pi_h \eta(z_h),\Pi_h \eta(z_h))  + b(\cdot,\Pi_h \eta(z_h),\Pi_h \eta(z_h)) \ds\\
 &\leq  \int_{t_0}^t  b(\cdot,\Pi_h \eta(z_h),\Pi_h \eta(z_h)) \ds \qquad \text{ for any }  t,t_0 \in I,
 \end{aligned}  
 \end{equation} 
 compare with Proposition~\ref{prop:discr-energy}. 
 This concludes the discussion on the level of discrete function spaces. 
 Below we also detail the structure on the level of coefficient vectors, since this is crucial for the implementation of our method.

 \subsubsection*{Details of the implementation}
 Both the system~\eqref{eq:pH-space-discrete} as well as the energy consistency can be formulated in the coefficient vectors $w_h$ with $z_h = I_h(w_h)$, or $c_h(z_h) = w_h$. 
 Provided that $j$, $r$ and $b$ are linear in the last argument, the resulting system has the form~\eqref{eq:pH-findim}.  
Indeed, under these conditions we may identify the operators $\overline J$, $\overline R \colon \setR^{\dim \BS_h} \to \setR^{\dim \BS_h}$ and $\overline B \colon I \times \setR^{\dim \BS_h} \times \setR^{\dim \BS_h} \to \setR^{\dim \BS_h}$ by 
\begin{align*}
\skp{\overline J(v_h)}{\varphi_h}_{\ell^2} &\coloneqq j(I_h (v_h), I_h (\varphi_h)), \quad \quad 
\skp{\overline R(v_h)}{\varphi_h}_{\ell^2} \coloneqq r(I_h (v_h), I_h (\varphi_h)), \quad \text{ and } \\
\skp{\overline B(t,v_h)}{\varphi_h}_{\ell^2} &\coloneqq b(t,I_h(v_h), I_h (\varphi_h)),
\end{align*}
for all $v_h,\varphi_h \in \setR^{\dim \BS_h}$. 
With $M$ the mass matrix of the basis of $\BS_h$ in the inner product of $Z$ the system reduces to 
\begin{subequations}\label{eq:pH-space-discrete-2}
\begin{align}
M 	\partial_t w_h  &= 
	\overline J(c_h(\Pi_h \eta(I_h (w_h)))) -\overline R(c_h(\Pi_h \eta(I_h (w_h))))
	 + \overline	B(\cdot, c_h(\Pi_h \eta(I_h (w_h)))),\\
	 w_h(0) &= c_h(\Pi_h z_0). 
	\end{align}
	\end{subequations}
	Noting that $M$ is regular, symmetric and positive definite, and setting $J(w_h) = M^{-1}\overline{J}(M^{-1} w_h)$, and similarly for the other operators, and choosing $\eta_h(w_h) = M c_h(\Pi_h \eta(I_h w_h))$ as above, system~\eqref{eq:pH-space-discrete-2} has the form~\eqref{eq:pH-findim} with Hamiltonian $\mathcal{H}_h(w_h) = \mathcal{H}(I_h w_h)$. 
	In particular, Assumption~\ref{as:problem-findim} is satisfied, and the discretization preserves the Hamiltonian structure.  

In all examples discussed below the Hilbert space $Z$ is $L^2(\Omega)$ or a vector-valued version thereof, and the forms $j,r,b$ are linear in the last argument. 
Note that we do not address well-posedness of the fully discrete problem here, since this would require more structure.

Similarly as before, we approximate the $L^\infty(I;\setR^\delta)$-norm of the error $w_{h} - w_{h\tau}$ with weighted Euclidean norm on $\setR^\delta$ (this is equivalent to the $L^{\infty}(I;L^2(\Omega))$-norm of $z_h - z_{h\tau}$), by evaluation on a time grid $\mathcal{T}_{\mathrm{fine}}  \coloneqq \{ i \tauref\colon i \in \{0, ..., 38913\}\}$ with time step size $\tauref \approx 1.29 \cdot 10^{-4}$, and additionally we consider the error at the time grid points. This means that we work with the relative errors 
	\begin{align}\label{def:error-Ewt}
		E \coloneqq \frac{\max_{t \in \mathcal{T}_{\mathrm{fine}}} \norm{w_{h\tau}(t) - w_h(t)}}{\max_{t \in \mathcal{T}_{\mathrm{fine}}} \norm{ w_h(t)}}  \quad \text{ and } \quad 
		E_{\tau} \coloneqq \frac{\max_{t \in \{t_0,t_1 , \ldots, t_m\}} \norm{w_{h\tau} (t)- w_h(t)}}{\max_{t \in \{t_0,t_1, \ldots, t_m\}} \norm{ w_h(t)}}. 
\end{align}

\subsection{Quasilinear wave equation}  
\noindent We consider Example~\ref{ex:wave-eq} on $Q = I \times \Omega$ for the one dimensional  domain $\Omega = [0,\ell] = [0,10]$ and for the time interval $I  = [0,T] = [0,5]$, for constants $\nu,\gamma \geq 0$, the pressure law $p(\rho) = \rho + \rho^3$, and the friction term
\begin{align}\label{eq:friction}
	F(v)  = \mathrm{sign}(v) \sqrt{|v|} = \abs{v}^{-1/2} v.
\end{align}
This corresponds to \eqref{eq:friction-powerlaw} with $s = 3/2$. 
System~\eqref{eq:wave-general} is supplemented with the boundary conditions 
\begin{equation*}
  p(\rho(t,0)) - \nu \del_x v(t,0) = g_0(t) \quad \text{and} \quad  p(\rho(t,\ell)) - \nu \del_x v(t,\ell) = g_\ell(t),
\end{equation*}
for some $g_0, g_\ell \in C(I)$, cf.~\eqref{eq:bc-wave-1}, \eqref{eq:bc-wave-2} with $\Gamma = \emptyset$. 
Recall also, that 
\begin{align*}
\eta(z) = \begin{pmatrix}
	p(\rho)\\ v 
\end{pmatrix} \quad \text{ for } z = \begin{pmatrix} \rho\\ v \end{pmatrix}. 
\end{align*}
	We employ the space discretization as described in Section~\ref{sec:space-discr} with piecewise  constant functions for $\rho$ and continuous piecewise linear functions for~$v$ to obtain a semi-discretization of~\eqref{eq:wave-general} in space. 
	We choose the Lagrange basis functions and the corresponding Lagrange interpolation. 
	When using $N$ equidistant inner grid points in $\Omega = [0,\ell]$
	with $h= \ell/(N+1)$ the numerical approximation $z_h = (\rho_h,v_h)\tp$ is represented by the coefficient vector $w_h = (w_{h,1}, w_{h,2})\tp \colon I \to \setR^{2N+3}$. 
Then, the method has the form~\eqref{eq:pH-space-discrete} with operators $j,r$ as specified in \eqref{eq:forms-wave-j}, \eqref{eq:forms-wave-r} and with $u = (g_0,g_{\ell})^{\tp}$ defining $b$ as in \eqref{eq:forms-wave-b} or \eqref{eq:forms-wave-b-v}, respectively. 
Note that we solve for the coefficients $w_h$ and hence we use the corresponding formulation of the method, cf.~\eqref{eq:pH-space-discrete-2}. 
In the following we use $N = 9$ interior discretization points, i.e., the space discretization is based on $10$ subintervals of equidistant length $h = 1$, and hence the system of form~\eqref{eq:pH-space-discrete} has dimension $\delta = 21$. 

For the time discretization for some $m \in \mathbb{N}$ we consider the time grid points 
$t_0,\dots,t_m$ with $\tau = T/m$ and $t_i = i \tau$ for any $i \in \{0, \ldots, m\}$. 
As parameters we choose $\gamma  = 0.1$ and $\nu \in \{0,1\}$. 
The approximate solutions $z_{h\tau}$ are obtained by applying Scheme~\ref{def:scheme} to the semi-discrete system of form~\eqref{eq:pH-space-discrete}, and analogously $w_{h\tau}$ results from applying Scheme~\ref{def:scheme} to \eqref{eq:pH-space-discrete-2}. 

With pressure law $p(\rho) = \rho + \rho^3$ it follows that $p(z_{h\tau}) \in \rV_{3k}^c(I_\tau;\setR^\delta)$ for  $z_{h\tau} \in \rV_{k}^c(I_\tau;\setR^{\delta})$. 
Thus, the integrals in the computation of $\proj\eta(z)$ are exact if~$s_\Pi = 2k$ nodes are used in the Gauß quadrature, which is used in the following.
Furthermore, as before Gauß quadrature with $s_Q = k$ nodes is used for the quadrature $Q_i$ in~\eqref{eq:scheme_local}. 

For the convergence analysis we use a manufactured solution. 
We choose the space interval midpoint values, and the space grid point values, respectively, of the functions 
\begin{equation*}
  \rho(t,x) \coloneqq v(t,x) \coloneqq \sin(t) \sin(x), 
\end{equation*}
as $w_h = (w_{h,1}, w_{h,2})\tp$. 
As defined in \eqref{def:error-Ewt} we consider the error terms comparing $w_h$ and $w_{h\tau}$. 
I.e.,  we only investigate the error due to the time discretization, and not the one arising by the space discretization. 
Figures~\ref{fig:damped_wave_nu0_irregular_friction_varying_degree} and \ref{fig:damped_wave_nu1_irregular_friction_varying_degree} show convergence of $E$ in $\tau$ for fixed space discretization with $h = 1$ for polynomial degrees~$k \in \{2,4,6\}$, for $s_Q = k$, for $\nu= 0$ and for $\nu = 1$, respectively. 
We observe optimal convergence order $\tau^{k+1}$ independently of the viscosity  $\nu=0$ and $\nu=1$. 
Also in the subsequent numerical experiments there is no difference between $\nu = 0$ and $\nu = 1$, which is why we omit the case $\nu=0$ in the following. 
Note that the manufactured solution satisfies that $v(t,\cdot) = 0$ for any $t \in \{\pi,2\pi, \ldots\}$. 
Seemingly this does not cause any difficulties with the friction term, which exemplifies that indeed our method can handle such situations, cf.~Remark~\ref{rmk:discr-comp}. 
In Figure~\ref{fig:damped_wave_nu1_irregular_friction_varying_degree_different_sampling} the same parameters and settings are used, to investigate nodal superconvergence. 
As before, we observe that the convergence rate is $\tau^{2k}$ in this case. 

In Figure~\ref{fig:damped_wave_nu1_irregular_friction_varying_discretization} convergence in $\tau$ is shown for $k = 4 = \quadnodes$, and $\projnodes = 2k$ and several spatial mesh sizes $h\in\{\tfrac{10}{9}, \tfrac{10}{17}, \tfrac{10}{33}, \tfrac{10}{65}\}$ for $\nu=1$. 
We observe that the convergence order in $\tau$ and also the error is independent of the spatial mesh size.  

To verify the energy consistency, in Figure~\ref{fig:damped_wave_nu1_irregular_friction_energybalance} the relative errors in the energy balance $\mathcal{E}$ as defined in \eqref{def:E} for the space discrete Hamiltonian $\mathcal{H}_h$ with $h = 1$ is depicted for $\nu=1$ and $\tau = 10^{-2}$.  
Here we use the boundary data 
\begin{equation*}
 g_0(t) = g_\ell(t) = 1 - \sin(t)
\end{equation*}
and as initial data for $w_h  = (w_{h,1},w_{h,2})\tp$ we use the corresponding point evaluations of the functions  
\begin{equation*}
  \rho(0,x) = 1 + \frac12 \sin\left(\frac{\pi x}{\ell}\right) ~\quad \text{and}\quad  v(0,x) = \left(\frac{4x}{\ell} - 2\right)^3,
\end{equation*}
as above. 
As before, we compare several polynomial degrees $k \in \{1,2,3,4\}$, with $s_Q = k$ and $s_\Pi = 2k$ and observe that our method satisfies the energy balance up to machine precision. 
This is expected because the integration to compute $\Pi$ is exact. 

\begin{figure}
	\hfill
	\begin{subfigure}[t]{.3\textwidth}
		\centering
		\scalebox{0.31}{\input{figures/damped_wave_nu0_irregular_friction_varying_degree.pgf}}
		\caption{
			convergence ($\nu=0$)}
		\label{fig:damped_wave_nu0_irregular_friction_varying_degree}
	\end{subfigure}
	\hfill 
	\begin{subfigure}[t]{.3\textwidth}
		\centering
		\scalebox{0.31}{\input{figures/damped_wave_nu1_irregular_friction_varying_degree.pgf}}
		\caption{
			convergence ($\nu=1$)}
		\label{fig:damped_wave_nu1_irregular_friction_varying_degree}
	\end{subfigure}
	\hfill 
	\begin{subfigure}[t]{.31\textwidth}
		\centering
		\scalebox{0.31}{\input{figures/damped_wave_nu1_irregular_friction_varying_degree_different_sampling.pgf}}
		\caption{
			nodal superconvergence ($\nu=1$)}
		\label{fig:damped_wave_nu1_irregular_friction_varying_degree_different_sampling}
	\end{subfigure}
	\hfill
	\caption{
		Convergence \eqref{fig:damped_wave_nu0_irregular_friction_varying_degree}  \eqref{fig:damped_wave_nu1_irregular_friction_varying_degree}
		and nodal superconvergence \eqref{fig:damped_wave_nu1_irregular_friction_varying_degree_different_sampling} in $\tau$
		for several polynomial degrees $k$,  $s_Q = k$ and $s_\Pi = 2k$ for the space-discrete quasilinear wave equation with fixed mesh size $h=1$ and $\nu\in\{0,1\}$.
	}
\end{figure}

\begin{figure}
    \centering
	\begin{subfigure}[t]{.3\textwidth}
		\centering
		\scalebox{0.31}{\input{figures/damped_wave_nu1_irregular_friction_varying_discretization.pgf}}
		\caption{
			convergence for $k = 4 = \quadnodes$, $\projnodes = 2k$, for several fixed $h$}
		\label{fig:damped_wave_nu1_irregular_friction_varying_discretization}
	\end{subfigure}
	\hspace{0.8em}
	\begin{subfigure}[t]{.3\textwidth}
		\centering
		\scalebox{0.31}{\input{figures/damped_wave_nu1_irregular_friction_energybalance.pgf}}
		\caption{
			relative energy error \eqref{def:E}, $\tau = 10^{-2}$, $h = 1$ for several $k = \quadnodes$ and $\projnodes = 2k$}
		\label{fig:damped_wave_nu1_irregular_friction_energybalance}
	\end{subfigure}
	\caption{
		Convergence in $\tau$ for several space-discretization parameters $h$ (\ref{fig:damped_wave_nu1_irregular_friction_varying_discretization})
		and 
		relative error in the energy balance with Hamiltonian $\mathcal{H}_h$  (\ref{fig:damped_wave_nu1_irregular_friction_energybalance})
		for the space-discrete quasilinear wave equation for $\nu=1$.
	}
\end{figure}

\subsection{Porous medium equation}
We consider the porous medium equation~\eqref{eq:parab-doubly-trafo} in Example~\ref{ex:doubly-parab} with $\beta(s) = s$ (corresponding to $p = 2$) and 
\begin{equation}\label{eq:porous-medium-aeps}
	\eta(z) = \alpha^{-1}(z) = |z|^{q-1} z + \eps z 
\end{equation}
for some $q>1$ and regularization parameter $\eps>0$. 
This results in the regularized porous medium equation 
\begin{align}\label{eq:porous-med}
		\partial_t z - \Delta(\alpha^{-1}(z)) 
		= \partial_t z - \div\left((q \abs{z}^{q-1} + \eps) \nabla z\right)  =  0.
\end{align}
We consider a one dimensional domain $\Omega=[0,\ell] = [0,15]$ and the time interval $I = [0,T] = [0,5]$. 

We employ the space discretization scheme as described in Section~\ref{sec:space-discr} with continuous piecewise linear functions for $z$.
When using $N$ equidistant inner grid points in $\Omega = [0,\ell]$ with $h = \ell/(N+1)$ the numerical approximation $z_h$ is represented by the coefficient vector $w_h\colon I \to \setR^{N+2}$.
In the following we use $N=14$ interior discretization points, i.e., 
we have $15$ subintervals of equidistant length $h=1$. 
This results in a system of form~\eqref{eq:pH-space-discrete} of dimension $\delta = 16$. 

For the time discretization of $I = [0,T]$ with time step size $\tau >0$ we consider the time grid points $t_i = i \tau$ for $i \in \{0, \ldots, m \}$ with $\tau = T/m$ for $m \in \mathbb{N}$. 

As parameters we choose $q \in \{1.5, 2, 3\}$ and $\eps \in \{10^{-10},10^{-8}\}$. 
The approximate solutions are obtained by applying Scheme~\ref{def:scheme}  to the semi-discrete system of form~\eqref{eq:pH-space-discrete}, or rather the corresponding system for $w_{h\tau}$. 

For $q=3$ the integrals in the computation of $\Pi \eta(z)$ are exact if $\projnodes = 2k$ nodes are used in the Gauß quadrature, which is used	in all experiments.
Furthermore, as before Gauß quadrature with $\quadnodes = k$ nodes is used for the quadrature $Q_i$ in~\eqref{eq:scheme_local}.

For the convergence analysis, we 
 work with a smooth manufactured solution $\zs$. 
Additionally we test our method for the Barenblatt solution $z_B$, which is an explicit non-smooth solution to equation~\eqref{eq:porous-med} without regularization, i.e., for $\eps = 0$.  
For the first test case we consider $z_h$ to be the continuous, piecewise affine interpolant of the smooth function 
\begin{equation*}
\zs(t,x) = \cos(t) \sin(x)
\end{equation*}
at the space grid points and we determine the corresponding right-hand side such that $z_h$ is the exact solution to the space-discrete regularized porous medium equation~\eqref{eq:porous-med} with $\eps >0$.
This means that for fixed regularization parameter $\eps>0$ we compare the space-discrete solution to the fully discrete solution, or rather we compare their coefficient vectors. 
For the latter test case, we consider the 1D Barenblatt solution, see, e.g.,~\cite[Sec.~1.2.2]{Vazquez2007}, which reads 
\begin{equation*}
	z_{\textup{B}}(t,x) = (t+1)^{-\frac{1}{q+1}} \max \left(0, 1 - \frac{q-1}{2q(q+1)} \frac{(x-\tfrac{\ell}{2})^2}{(t+1)^{\frac{2}{q+1}}} \right)^{\frac{1}{q-1}}.
\end{equation*}
In this case $z_h$ is the continuous piecewise affine interpolant of $z_B$, and $w_h$ is the corresponding coefficient vector. 
We compare this interpolation of the Barenblatt solution (which can be expected to be close to the solution to the space-discrete equation for $\eps = 0$) and the numerical solution to the fully-discrete method for the regularized porous medium equation depending on $\eps$. 
Hence, both the time-discretization and the regularization as well as the (fixed) space discretization affect the error.
As before, we compute the errors $E$ and $E_\tau$, as in \eqref{def:error-Ewt}. 

\begin{table}
	\hfill
	\begin{subtable}[t]{.49\textwidth}
		\centering
		\scalebox{0.8}{
\begin{tabular}{|c|cc|cc|}
\hline
    \multirow{2}{\widthof{$\tau$}}{$\tau$} & \multicolumn{2}{|c}{$k=2$} & \multicolumn{2}{|c|}{$k=4$}\Tstrut \\
     & error $E_{}$ & EOC & error $E_{}$ & EOC\Bstrut \\ \hline
    $2.50 \cdot 10^{-01}$ & $3.70 \cdot 10^{-04}$ & - & $1.19 \cdot 10^{-04}$ & -\Tstrut \\
    $1.28 \cdot 10^{-01}$ & $1.72 \cdot 10^{-04}$ & $1.15$ & $3.70 \cdot 10^{-05}$ & $1.75$ \\
    $6.49 \cdot 10^{-02}$ & $1.22 \cdot 10^{-05}$ & $3.89$ & $6.90 \cdot 10^{-06}$ & $2.47$ \\
    $3.27 \cdot 10^{-02}$ & $6.29 \cdot 10^{-06}$ & $0.97$ & $1.67 \cdot 10^{-06}$ & $2.06$ \\
    $1.64 \cdot 10^{-02}$ & $2.45 \cdot 10^{-06}$ & $1.37$ & $1.17 \cdot 10^{-06}$ & $0.52$ \\
    $8.21 \cdot 10^{-03}$ & $4.71 \cdot 10^{-07}$ & $2.38$ & $2.04 \cdot 10^{-07}$ & $2.53$ \\
    $4.11 \cdot 10^{-03}$ & $1.27 \cdot 10^{-07}$ & $1.89$ & $4.55 \cdot 10^{-08}$ & $2.16$ \\
    $2.06 \cdot 10^{-03}$ & $3.77 \cdot 10^{-08}$ & $1.76$ & $6.80 \cdot 10^{-09}$ & $2.74$ \\
    $1.03 \cdot 10^{-03}$ & $8.71 \cdot 10^{-09}$ & $2.11$ & $1.19 \cdot 10^{-09}$ & $2.51$\Bstrut \\
\hline
      & & $\hspace{-11pt}\diameter\; 1.94$ & & $\hspace{-11pt}\diameter\; 2.09$\Tstrut\Bstrut \\ \hline
\end{tabular}
}
		\caption{
			Comparison of space-interpolant of $z_{\text{B}}$ ($\varepsilon = 0$) \\and $\varepsilon$-dependent fully discrete solution}
			\label{fig:porous_medium_q1.5_eps1e-10_varying_degree-zB}
	\end{subtable}
	\hfill
	\begin{subtable}[t]{.49\textwidth}
		\centering
		\scalebox{0.8}{
\begin{tabular}{|c|cc|cc|}
\hline
    \multirow{2}{\widthof{$\tau$}}{$\tau$} & \multicolumn{2}{|c}{$k=2$} & \multicolumn{2}{|c|}{$k=4$}\Tstrut \\
     & error $E_{}$ & EOC & error $E_{}$ & EOC\Bstrut \\ \hline
    $2.50 \cdot 10^{-01}$ & $3.36 \cdot 10^{-04}$ & - & $8.51 \cdot 10^{-05}$ & -\Tstrut \\
    $1.28 \cdot 10^{-01}$ & $6.23 \cdot 10^{-05}$ & $2.52$ & $1.63 \cdot 10^{-05}$ & $2.47$ \\
    $6.49 \cdot 10^{-02}$ & $2.40 \cdot 10^{-05}$ & $1.41$ & $1.09 \cdot 10^{-05}$ & $0.59$ \\
    $3.27 \cdot 10^{-02}$ & $7.27 \cdot 10^{-06}$ & $1.74$ & $1.54 \cdot 10^{-06}$ & $2.85$ \\
    $1.64 \cdot 10^{-02}$ & $7.06 \cdot 10^{-07}$ & $3.38$ & $1.16 \cdot 10^{-07}$ & $3.75$ \\
    $8.21 \cdot 10^{-03}$ & $3.22 \cdot 10^{-08}$ & $4.47$ & $3.41 \cdot 10^{-09}$ & $5.10$ \\
    $4.11 \cdot 10^{-03}$ & $7.92 \cdot 10^{-09}$ & $2.02$ & $5.98 \cdot 10^{-10}$ & $2.52$ \\
    $2.06 \cdot 10^{-03}$ & $2.03 \cdot 10^{-09}$ & $1.96$ & $4.78 \cdot 10^{-10}$ & $0.32$ \\
    $1.03 \cdot 10^{-03}$ & $7.31 \cdot 10^{-10}$ & $1.48$ & $1.30 \cdot 10^{-10}$ & $1.88$\Bstrut \\
\hline
      & & $\hspace{-11pt}\diameter\; 2.37$ & & $\hspace{-11pt}\diameter\; 2.44$\Tstrut\Bstrut \\ \hline
\end{tabular}
}
			\caption{
				Comparison of space-interpolant of $\zs$ and  fully \\
			discrete solution (both $\varepsilon$-dependent)}
			\label{fig:porous_medium_q1.5_eps1e-10_varying_degree-zs}
	\end{subtable}
	\hfill
	\caption{
		Convergence in $\tau$ for several polynomial degrees $k$, 
		for $\quadnodes = k$, and $\projnodes = 2k$, for the space-discrete (fixed mesh size $h = 1$) regularized ($\eps=10^{-10}$) porous medium equation for $q=1.5$.	
	}
	\label{fig:porous_medium_q1.5_eps1e-10_varying_degree}
\end{table}
\begin{table}[h]
\hfill
	\begin{subtable}[t]{.49\textwidth}
		\centering
		\scalebox{0.8}{
\begin{tabular}{|c|cc|cc|}
\hline
    \multirow{2}{\widthof{$\tau$}}{$\tau$} & \multicolumn{2}{|c}{$k=2$} & \multicolumn{2}{|c|}{$k=4$}\Tstrut \\
     & error $E_{}$ & EOC & error $E_{}$ & EOC\Bstrut \\ \hline
    $2.50 \cdot 10^{-01}$ & $3.70 \cdot 10^{-04}$ & - & $1.19 \cdot 10^{-04}$ & -\Tstrut \\
    $1.28 \cdot 10^{-01}$ & $1.72 \cdot 10^{-04}$ & $1.15$ & $3.70 \cdot 10^{-05}$ & $1.75$ \\
    $6.49 \cdot 10^{-02}$ & $1.22 \cdot 10^{-05}$ & $3.89$ & $6.90 \cdot 10^{-06}$ & $2.47$ \\
    $3.27 \cdot 10^{-02}$ & $6.29 \cdot 10^{-06}$ & $0.97$ & $1.67 \cdot 10^{-06}$ & $2.06$ \\
    $1.64 \cdot 10^{-02}$ & $2.45 \cdot 10^{-06}$ & $1.37$ & $1.17 \cdot 10^{-06}$ & $0.52$ \\
    $8.21 \cdot 10^{-03}$ & $4.71 \cdot 10^{-07}$ & $2.38$ & $2.04 \cdot 10^{-07}$ & $2.53$ \\
    $4.11 \cdot 10^{-03}$ & $1.27 \cdot 10^{-07}$ & $1.89$ & $4.55 \cdot 10^{-08}$ & $2.16$ \\
    $2.06 \cdot 10^{-03}$ & $3.77 \cdot 10^{-08}$ & $1.76$ & $6.80 \cdot 10^{-09}$ & $2.74$ \\
    $1.03 \cdot 10^{-03}$ & $8.71 \cdot 10^{-09}$ & $2.11$ & $1.19 \cdot 10^{-09}$ & $2.51$\Bstrut \\
\hline
      & & $\hspace{-11pt}\diameter\; 1.94$ & & $\hspace{-11pt}\diameter\; 2.09$\Tstrut\Bstrut \\ \hline
\end{tabular}
}
		\caption{$z_{\text{B}}$}
	\end{subtable}
	\hfill
	\begin{subtable}[t]{.49\textwidth}
		\centering
		\scalebox{0.8}{
\begin{tabular}{|c|cc|cc|}
\hline
    \multirow{2}{\widthof{$\tau$}}{$\tau$} & \multicolumn{2}{|c}{$k=2$} & \multicolumn{2}{|c|}{$k=4$}\Tstrut \\
     & error $E_{}$ & EOC & error $E_{}$ & EOC\Bstrut \\ \hline
    $2.50 \cdot 10^{-01}$ & $3.36 \cdot 10^{-04}$ & - & $8.51 \cdot 10^{-05}$ & -\Tstrut \\
    $1.28 \cdot 10^{-01}$ & $6.23 \cdot 10^{-05}$ & $2.52$ & $1.63 \cdot 10^{-05}$ & $2.47$ \\
    $6.49 \cdot 10^{-02}$ & $2.40 \cdot 10^{-05}$ & $1.41$ & $1.09 \cdot 10^{-05}$ & $0.59$ \\
    $3.27 \cdot 10^{-02}$ & $7.27 \cdot 10^{-06}$ & $1.74$ & $1.54 \cdot 10^{-06}$ & $2.85$ \\
    $1.64 \cdot 10^{-02}$ & $7.06 \cdot 10^{-07}$ & $3.38$ & $1.16 \cdot 10^{-07}$ & $3.75$ \\
    $8.21 \cdot 10^{-03}$ & $3.22 \cdot 10^{-08}$ & $4.47$ & $3.41 \cdot 10^{-09}$ & $5.10$ \\
    $4.11 \cdot 10^{-03}$ & $7.92 \cdot 10^{-09}$ & $2.02$ & $5.98 \cdot 10^{-10}$ & $2.52$ \\
    $2.06 \cdot 10^{-03}$ & $2.03 \cdot 10^{-09}$ & $1.96$ & $4.78 \cdot 10^{-10}$ & $0.32$ \\
    $1.03 \cdot 10^{-03}$ & $7.31 \cdot 10^{-10}$ & $1.48$ & $1.30 \cdot 10^{-10}$ & $1.88$\Bstrut \\
\hline
      & & $\hspace{-11pt}\diameter\; 2.37$ & & $\hspace{-11pt}\diameter\; 2.44$\Tstrut\Bstrut \\ \hline
\end{tabular}
}
		\caption{$\zs$}
	\end{subtable}
	\hfill
	\caption{
		Convergence in $\tau$ for several polynomial degrees $k$, 
		for $\quadnodes = k$, and $\projnodes = 2k$, for the space-discrete (fixed mesh size $h = 1$) regularized ($\eps=10^{-8}$) porous medium equation for $q=1.5$.
	}
	\label{fig:porous_medium_q1.5_eps1e-08_varying_degree}
\end{table}

To investigate the convergence order in $\tau$ for fixed space discretization, we compute the estimated order of convergence (eoc) of the errors. 
For a sequence of $(\tau_j)_j$ and errors $(e_j)_j$ it is defined by
\begin{align}
	\texttt{eoc}_j \coloneqq \frac{\log\left(\frac{e_j}{e_{j-1}}\right)}{\log\left(\frac{\tau_{j}}{ \tau_{j-1}}\right)}, \qquad \text{ for } j \in \mathbb{N}. 
\end{align}
We compute the eoc for errors above the threshold $10^{-14}$ and of all those we also compute the arithmetic mean, see Tables~\ref{fig:porous_medium_q1.5_eps1e-10_varying_degree}--\ref{fig:porous_medium_q3.0_eps1e-10_varying_degree_different_sampling}. 

Tables~\ref{fig:porous_medium_q1.5_eps1e-10_varying_degree} and~\ref{fig:porous_medium_q1.5_eps1e-08_varying_degree} show the error $E$ for $q=1.5$, for polynomial degrees~$k \in \{2,4\}$ and for  $\quadnodes = k$, $\projnodes = 2k$, for the solutions $z_{\text{B}}$ and $\zs$, each, and for regularization parameters $\eps = 10^{-10}$ and $\eps =  10^{-8}$, respectively. 
Evidently, the choice of the regularization parameter does not have any effect, since both tables show the same results. 
For this reason we shall refrain from presenting the results for both choices in the following, and we only include the results for $\eps = 10^{-10}$. 
The same parameters and settings ($\eps=10^{-10}$) are used for $q = 2$ in Table~\ref{fig:porous_medium_q2.0_eps1e-10_varying_degree} and for $q = 3$ in Table~\ref{fig:porous_medium_q3.0_eps1e-10_varying_degree}. 
We observe for the Barenblatt solution $z_B$ that the convergence order in $\tau$ is reduced for both $k = 2$ and $k = 4$, and it is the  smaller, the larger $q$ is, see the average values in~Tables~\ref{fig:porous_medium_q1.5_eps1e-10_varying_degree-zB}, 			\ref{fig:porous_medium_q2.0_eps1e-10_varying_degree-zB} and \ref{fig:porous_medium_q3.0_eps1e-10_varying_degree-zB}. 
This can indeed be expected due to the lack of higher regularity of the Barenblatt solution. 
More precisely, solutions to the porous medium equation have less regularity in time the larger $q$ is, see \cite{GessSauerTadmor2020}. 
On the other hand, for the smooth solution~$\zs$ higher convergence order can be observed, see Tables~\ref{fig:porous_medium_q1.5_eps1e-10_varying_degree-zs}, 	\ref{fig:porous_medium_q2.0_eps1e-10_varying_degree-zs} and 			\ref{fig:porous_medium_q3.0_eps1e-10_varying_degree-zs}. 
However, the optimal order $\tau^{k+1}$ is only obtained for $q = 3$. 
This can be expected since $\eta$ is smooth for $q = 3$, but only $C^1$ for $q \in \{1.5,2\}$.  


\begin{table}
\hfill
	\begin{subtable}[t]{.49\textwidth}
		\centering
		\scalebox{0.8}{
\begin{tabular}{|c|cc|cc|}
\hline
    \multirow{2}{\widthof{$\tau$}}{$\tau$} & \multicolumn{2}{|c}{$k=2$} & \multicolumn{2}{|c|}{$k=4$}\Tstrut \\
     & error $E_{}$ & EOC & error $E_{}$ & EOC\Bstrut \\ \hline
    $2.50 \cdot 10^{-01}$ & $1.48 \cdot 10^{-02}$ & - & $8.52 \cdot 10^{-03}$ & -\Tstrut \\
    $1.28 \cdot 10^{-01}$ & $1.79 \cdot 10^{-02}$ & - & $4.35 \cdot 10^{-03}$ & $1.01$ \\
    $6.49 \cdot 10^{-02}$ & $2.00 \cdot 10^{-03}$ & $3.22$ & $1.37 \cdot 10^{-03}$ & $1.69$ \\
    $3.27 \cdot 10^{-02}$ & $2.17 \cdot 10^{-03}$ & - & $6.96 \cdot 10^{-04}$ & $0.99$ \\
    $1.64 \cdot 10^{-02}$ & $6.96 \cdot 10^{-04}$ & $1.65$ & $6.48 \cdot 10^{-04}$ & $0.10$ \\
    $8.21 \cdot 10^{-03}$ & $6.96 \cdot 10^{-04}$ & - & $1.14 \cdot 10^{-04}$ & $2.52$ \\
    $4.11 \cdot 10^{-03}$ & $2.70 \cdot 10^{-04}$ & $1.37$ & $2.87 \cdot 10^{-04}$ & - \\
    $2.06 \cdot 10^{-03}$ & $2.51 \cdot 10^{-04}$ & $0.10$ & $1.08 \cdot 10^{-04}$ & $1.41$ \\
    $1.03 \cdot 10^{-03}$ & $5.63 \cdot 10^{-05}$ & $2.16$ & $3.04 \cdot 10^{-05}$ & $1.83$\Bstrut \\
\hline
      & & $\hspace{-11pt}\diameter\; 1.01$ & & $\hspace{-11pt}\diameter\; 1.03$\Tstrut\Bstrut \\ \hline
\end{tabular}
}
		\caption{$z_{\text{B}}$}
			\label{fig:porous_medium_q2.0_eps1e-10_varying_degree-zB}
	\end{subtable}
	\hfill
	\begin{subtable}[t]{.49\textwidth}
		\centering
		\scalebox{0.8}{
\begin{tabular}{|c|cc|cc|}
\hline
    \multirow{2}{\widthof{$\tau$}}{$\tau$} & \multicolumn{2}{|c}{$k=2$} & \multicolumn{2}{|c|}{$k=4$}\Tstrut \\
     & error $E_{}$ & EOC & error $E_{}$ & EOC\Bstrut \\ \hline
    $2.50 \cdot 10^{-01}$ & $1.51 \cdot 10^{-04}$ & - & $2.03 \cdot 10^{-05}$ & -\Tstrut \\
    $1.28 \cdot 10^{-01}$ & $3.94 \cdot 10^{-05}$ & $2.01$ & $4.94 \cdot 10^{-06}$ & $2.12$ \\
    $6.49 \cdot 10^{-02}$ & $7.25 \cdot 10^{-06}$ & $2.49$ & $1.85 \cdot 10^{-06}$ & $1.44$ \\
    $3.27 \cdot 10^{-02}$ & $1.44 \cdot 10^{-06}$ & $2.36$ & $1.74 \cdot 10^{-07}$ & $3.45$ \\
    $1.64 \cdot 10^{-02}$ & $1.48 \cdot 10^{-07}$ & $3.29$ & $1.40 \cdot 10^{-08}$ & $3.65$ \\
    $8.21 \cdot 10^{-03}$ & $4.46 \cdot 10^{-09}$ & $5.07$ & $1.02 \cdot 10^{-11}$ & $10.45$ \\
    $4.11 \cdot 10^{-03}$ & $5.57 \cdot 10^{-10}$ & $3.00$ & $3.16 \cdot 10^{-11}$ & - \\
    $2.06 \cdot 10^{-03}$ & $7.33 \cdot 10^{-11}$ & $2.93$ & $1.39 \cdot 10^{-11}$ & $1.19$ \\
    $1.03 \cdot 10^{-03}$ & $2.95 \cdot 10^{-11}$ & $1.31$ & $3.40 \cdot 10^{-12}$ & $2.03$\Bstrut \\
\hline
      & & $\hspace{-11pt}\diameter\; 2.81$ & & $\hspace{-11pt}\diameter\; 2.84$\Tstrut\Bstrut \\ \hline
\end{tabular}
}
		\caption{$\zs$}
			\label{fig:porous_medium_q2.0_eps1e-10_varying_degree-zs}
	\end{subtable}
	\hfill
	\caption{
		Convergence in $\tau$ for several polynomial degrees $k$, 
		for $\quadnodes = k$, and $\projnodes = 2k$, for the space-discrete (fixed mesh size $h = 1$) regularized ($\eps=10^{-10}$) porous medium equation for $q=2$. 		
	}
	\label{fig:porous_medium_q2.0_eps1e-10_varying_degree}
\end{table}

\begin{table}
	\hfill
	\begin{subtable}[t]{.49\textwidth}
		\centering
		\scalebox{0.8}{
\begin{tabular}{|c|cc|cc|}
\hline
    \multirow{2}{\widthof{$\tau$}}{$\tau$} & \multicolumn{2}{|c}{$k=2$} & \multicolumn{2}{|c|}{$k=4$}\Tstrut \\
     & error $E_{}$ & EOC & error $E_{}$ & EOC\Bstrut \\ \hline
    $2.50 \cdot 10^{-01}$ & $1.17 \cdot 10^{-01}$ & - & $5.37 \cdot 10^{-02}$ & -\Tstrut \\
    $1.28 \cdot 10^{-01}$ & $8.66 \cdot 10^{-02}$ & $0.44$ & $2.09 \cdot 10^{-01}$ & - \\
    $6.49 \cdot 10^{-02}$ & $1.78 \cdot 10^{-02}$ & $2.33$ & $6.50 \cdot 10^{-02}$ & $1.72$ \\
    $3.27 \cdot 10^{-02}$ & $4.42 \cdot 10^{-02}$ & - & $2.37 \cdot 10^{-02}$ & $1.47$ \\
    $1.64 \cdot 10^{-02}$ & $2.14 \cdot 10^{-02}$ & $1.05$ & $1.38 \cdot 10^{-02}$ & $0.78$ \\
    $8.21 \cdot 10^{-03}$ & $1.19 \cdot 10^{-02}$ & $0.84$ & $1.89 \cdot 10^{-02}$ & - \\
    $4.11 \cdot 10^{-03}$ & $2.20 \cdot 10^{-02}$ & - & $5.36 \cdot 10^{-03}$ & $1.82$ \\
    $2.06 \cdot 10^{-03}$ & $6.86 \cdot 10^{-03}$ & $1.68$ & $8.24 \cdot 10^{-03}$ & - \\
    $1.03 \cdot 10^{-03}$ & $6.51 \cdot 10^{-03}$ & $0.08$ & $5.62 \cdot 10^{-03}$ & $0.55$\Bstrut \\
\hline
      & & $\hspace{-11pt}\diameter\; 0.53$ & & $\hspace{-11pt}\diameter\; 0.40$\Tstrut\Bstrut \\ \hline
\end{tabular}
}
		\caption{$z_{\text{B}}$}
			\label{fig:porous_medium_q3.0_eps1e-10_varying_degree-zB}
	\end{subtable}
	\hfill
	\begin{subtable}[t]{.49\textwidth}
		\centering
		\scalebox{0.8}{
\begin{tabular}{|c|cc|cc|}
\hline
    \multirow{2}{\widthof{$\tau$}}{$\tau$} & \multicolumn{2}{|c}{$k=2$} & \multicolumn{2}{|c|}{$k=4$}\Tstrut \\
     & error $E_{}$ & EOC & error $E_{}$ & EOC\Bstrut \\ \hline
    $2.50 \cdot 10^{-01}$ & $1.89 \cdot 10^{-04}$ & - & $1.26 \cdot 10^{-07}$ & -\Tstrut \\
    $1.28 \cdot 10^{-01}$ & $2.01 \cdot 10^{-05}$ & $3.36$ & $2.06 \cdot 10^{-09}$ & $6.16$ \\
    $6.49 \cdot 10^{-02}$ & $2.37 \cdot 10^{-06}$ & $3.14$ & $3.72 \cdot 10^{-11}$ & $5.90$ \\
    $3.27 \cdot 10^{-02}$ & $2.90 \cdot 10^{-07}$ & $3.06$ & $1.00 \cdot 10^{-12}$ & $5.27$ \\
    $1.64 \cdot 10^{-02}$ & $3.60 \cdot 10^{-08}$ & $3.03$ & $3.06 \cdot 10^{-14}$ & $5.06$ \\
    $8.21 \cdot 10^{-03}$ & $4.47 \cdot 10^{-09}$ & $3.01$ & $2.62 \cdot 10^{-15}$ & $3.55$ \\
    $4.11 \cdot 10^{-03}$ & $5.58 \cdot 10^{-10}$ & $3.01$ & $2.49 \cdot 10^{-15}$ & - \\
    $2.06 \cdot 10^{-03}$ & $6.90 \cdot 10^{-11}$ & $3.02$ & $2.91 \cdot 10^{-15}$ & - \\
    $1.03 \cdot 10^{-03}$ & $8.49 \cdot 10^{-12}$ & $3.02$ & $5.10 \cdot 10^{-15}$ & -\Bstrut \\
\hline
      & & $\hspace{-11pt}\diameter\; 3.08$ & & $\hspace{-11pt}\diameter\; 5.19$\Tstrut\Bstrut \\ \hline
\end{tabular}
}
		\caption{$\zs$}
			\label{fig:porous_medium_q3.0_eps1e-10_varying_degree-zs}
	\end{subtable}
	\hfill
	\caption{
		Convergence in $\tau$ for several polynomial degrees $k$, 	for $\quadnodes = k$, and $\projnodes = 2k$, for the space-discrete (fixed mesh size $h = 1$) regularized ($\eps=10^{-10}$) porous medium equation for $q=3$. 
	}
	\label{fig:porous_medium_q3.0_eps1e-10_varying_degree}
\end{table}


In Tables~\ref{fig:porous_medium_q1.5_eps1e-10_varying_degree_different_sampling}--\ref{fig:porous_medium_q3.0_eps1e-10_varying_degree_different_sampling} we investigate nodal superconvergence for $q \in \{1.5, 2, 3\}$ and $\eps = 10^{-10}$.  
As before, we observe that the convergence rate is $\tau^{2k}$ only for $q=3$ for the smooth solution $\zs$ and that the convergence order is smaller for all other cases.

In Figure~\ref{fig:porous_medium_barenblatt_varying_discretization} and Figure~\ref{fig:porous_medium_smooth_varying_discretization} convergence is shown for several spatial mesh sizes 
$h\in\{\tfrac{15}{9}, \tfrac{15}{17}, \tfrac{15}{33}, \tfrac{15}{65}\}$, for $z_{\text{B}}$ and $\zs$, respectively.
Here we have chosen $k = 4 = \quadnodes$, $\projnodes = 2k = 8$,  $\eps=10^{-10}$ and all cases $q\in\{1.5, 2, 3\}$ are presented. 
We observe that the error does not increase when the number of inner grid points is increased.

To verify the energy consistency,  in Figure~\ref{fig:porous_medium_barenblatt_energybalance} and \ref{fig:porous_medium_smooth_energybalance} the relative errors in the energy balance as defined in~\eqref{def:E} for $\mathcal{H}_h$ for discretization parameter $h = 1$, for $q\in\{1.5, 2, 3\}$ and regularization parameter $\eps=10^{-10}$ and time step $\tau = 10^{-2}$ is presented. 
Here we use point evaluations of the function $z_{\text{B}}$ and $\zs$, respectively, as initial data for $w_h$. 
As before, we compare several polynomial degrees $k \in \{1,2,3,4\}$, and we choose $\quadnodes = k$ and $\projnodes = 2k$. 
We observe that for both the smooth and the non-smooth initial datum the energy balance is satisfied up to machine precision for the case where $\eta$ is smooth ($q = 3$), see Figures~\ref{fig:porous_medium_smooth_energybalance-q3} and~\ref{fig:porous_medium_barenblatt_energybalance-q3}, respectively. 
In the case of less regular $\eta$ the relative error in the energy balance is still mostly close to machine precision, but slightly larger at certain time points, especially for low polynomial degree $k$.
Since the energy balance is exact for exact integration of the $L^2$-projection, the reason for this is the lack of exactness of the quadrature to compute $\Pi$. 
To achieve a lower error in the energy balance it is inevitable to choose a 
 larger number of quadrature points $\projnodes$. 

\subsection*{Acknowledgements}
\vbox{
	All authors thank the Deutsche Forschungsgemeinschaft for their support within the subprojects B03, C05, and C09 in the Son\-der\-for\-schungs\-be\-reich/Trans\-re\-gio 154 ``Mathematical Modelling, Simulation and Optimization using the Example of Gas Networks'' (Project 239904186).
	Further, A.\,Karsai thanks P.\,Schulze for many helpful comments and discussions.
}

\begin{table}
\hfill
	\begin{subtable}[t]{.49\textwidth}
		\centering
		\scalebox{0.8}{
\begin{tabular}{|c|cc|cc|}
\hline
    \multirow{2}{\widthof{$\tau$}}{$\tau$} & \multicolumn{2}{|c}{$k=2$} & \multicolumn{2}{|c|}{$k=4$}\Tstrut \\
     & error $E_{\tau}$ & EOC & error $E_{\tau}$ & EOC\Bstrut \\ \hline
    $2.50 \cdot 10^{-01}$ & $1.46 \cdot 10^{-04}$ & - & $1.01 \cdot 10^{-04}$ & -\Tstrut \\
    $1.28 \cdot 10^{-01}$ & $1.25 \cdot 10^{-04}$ & $0.23$ & $3.05 \cdot 10^{-05}$ & $1.79$ \\
    $6.49 \cdot 10^{-02}$ & $5.86 \cdot 10^{-06}$ & $4.50$ & $5.92 \cdot 10^{-06}$ & $2.41$ \\
    $3.27 \cdot 10^{-02}$ & $5.87 \cdot 10^{-06}$ & - & $1.33 \cdot 10^{-06}$ & $2.17$ \\
    $1.64 \cdot 10^{-02}$ & $4.65 \cdot 10^{-07}$ & $3.68$ & $9.35 \cdot 10^{-07}$ & $0.52$ \\
    $8.21 \cdot 10^{-03}$ & $2.26 \cdot 10^{-07}$ & $1.04$ & $1.66 \cdot 10^{-07}$ & $2.50$ \\
    $4.11 \cdot 10^{-03}$ & $2.70 \cdot 10^{-08}$ & $3.07$ & $3.42 \cdot 10^{-08}$ & $2.28$ \\
    $2.06 \cdot 10^{-03}$ & $2.52 \cdot 10^{-08}$ & $0.10$ & $2.96 \cdot 10^{-09}$ & $3.53$ \\
    $1.03 \cdot 10^{-03}$ & $8.70 \cdot 10^{-09}$ & $1.54$ & $1.19 \cdot 10^{-09}$ & $1.31$\Bstrut \\
\hline
      & & $\hspace{-11pt}\diameter\; 1.77$ & & $\hspace{-11pt}\diameter\; 2.06$\Tstrut\Bstrut \\ \hline
\end{tabular}
}
		\caption{$z_{\text{B}}$}
	\end{subtable}
	\hfill
	\begin{subtable}[t]{.49\textwidth}
		\centering
		\scalebox{0.8}{
\begin{tabular}{|c|cc|cc|}
\hline
    \multirow{2}{\widthof{$\tau$}}{$\tau$} & \multicolumn{2}{|c}{$k=2$} & \multicolumn{2}{|c|}{$k=4$}\Tstrut \\
     & error $E_{\tau}$ & EOC & error $E_{\tau}$ & EOC\Bstrut \\ \hline
    $2.50 \cdot 10^{-01}$ & $3.36 \cdot 10^{-04}$ & - & $8.50 \cdot 10^{-05}$ & -\Tstrut \\
    $1.28 \cdot 10^{-01}$ & $4.59 \cdot 10^{-05}$ & $2.98$ & $7.21 \cdot 10^{-06}$ & $3.70$ \\
    $6.49 \cdot 10^{-02}$ & $2.33 \cdot 10^{-05}$ & $1.00$ & $8.62 \cdot 10^{-06}$ & - \\
    $3.27 \cdot 10^{-02}$ & $6.91 \cdot 10^{-06}$ & $1.77$ & $1.13 \cdot 10^{-06}$ & $2.96$ \\
    $1.64 \cdot 10^{-02}$ & $6.80 \cdot 10^{-08}$ & $6.70$ & $2.81 \cdot 10^{-08}$ & $5.35$ \\
    $8.21 \cdot 10^{-03}$ & $1.85 \cdot 10^{-08}$ & $1.88$ & $1.54 \cdot 10^{-09}$ & $4.20$ \\
    $4.11 \cdot 10^{-03}$ & $3.06 \cdot 10^{-09}$ & $2.60$ & $3.75 \cdot 10^{-10}$ & $2.04$ \\
    $2.06 \cdot 10^{-03}$ & $9.60 \cdot 10^{-10}$ & $1.68$ & $3.67 \cdot 10^{-10}$ & $0.03$ \\
    $1.03 \cdot 10^{-03}$ & $7.31 \cdot 10^{-10}$ & $0.39$ & $8.89 \cdot 10^{-11}$ & $2.05$\Bstrut \\
\hline
      & & $\hspace{-11pt}\diameter\; 2.37$ & & $\hspace{-11pt}\diameter\; 2.51$\Tstrut\Bstrut \\ \hline
\end{tabular}
}
		\caption{$\zs$}
	\end{subtable}
	\hfill
	\caption{
	Nodal superconvergence in $\tau$ for several polynomial degrees $k$, for $\quadnodes = k$, and $\projnodes = 2k$, for the space-discrete (fixed mesh size $h = 1$) regularized ($\eps=10^{-10}$) porous medium equation for $q=1.5$. 
	}
	\label{fig:porous_medium_q1.5_eps1e-10_varying_degree_different_sampling}
\end{table}

\begin{table}
	\hfill
	\begin{subtable}[t]{.49\textwidth}
		\centering
		\scalebox{0.8}{
\begin{tabular}{|c|cc|cc|}
\hline
    \multirow{2}{\widthof{$\tau$}}{$\tau$} & \multicolumn{2}{|c}{$k=2$} & \multicolumn{2}{|c|}{$k=4$}\Tstrut \\
     & error $E_{\tau}$ & EOC & error $E_{\tau}$ & EOC\Bstrut \\ \hline
    $2.50 \cdot 10^{-01}$ & $1.29 \cdot 10^{-02}$ & - & $5.92 \cdot 10^{-03}$ & -\Tstrut \\
    $1.28 \cdot 10^{-01}$ & $1.53 \cdot 10^{-02}$ & - & $4.28 \cdot 10^{-03}$ & $0.48$ \\
    $6.49 \cdot 10^{-02}$ & $1.37 \cdot 10^{-03}$ & $3.55$ & $1.13 \cdot 10^{-03}$ & $1.96$ \\
    $3.27 \cdot 10^{-02}$ & $2.03 \cdot 10^{-03}$ & - & $6.96 \cdot 10^{-04}$ & $0.70$ \\
    $1.64 \cdot 10^{-02}$ & $6.96 \cdot 10^{-04}$ & $1.55$ & $6.48 \cdot 10^{-04}$ & $0.10$ \\
    $8.21 \cdot 10^{-03}$ & $6.96 \cdot 10^{-04}$ & - & $1.05 \cdot 10^{-04}$ & $2.63$ \\
    $4.11 \cdot 10^{-03}$ & $2.70 \cdot 10^{-04}$ & $1.37$ & $2.72 \cdot 10^{-04}$ & - \\
    $2.06 \cdot 10^{-03}$ & $2.14 \cdot 10^{-04}$ & $0.33$ & $1.02 \cdot 10^{-04}$ & $1.43$ \\
    $1.03 \cdot 10^{-03}$ & $5.63 \cdot 10^{-05}$ & $1.93$ & $2.79 \cdot 10^{-05}$ & $1.86$\Bstrut \\
\hline
      & & $\hspace{-11pt}\diameter\; 0.99$ & & $\hspace{-11pt}\diameter\; 0.97$\Tstrut\Bstrut \\ \hline
\end{tabular}
}
		\caption{$z_{\text{B}}$}
	\end{subtable}
\hfill
	\begin{subtable}[t]{.49\textwidth}
		\centering
		\scalebox{0.8}{
\begin{tabular}{|c|cc|cc|}
\hline
    \multirow{2}{\widthof{$\tau$}}{$\tau$} & \multicolumn{2}{|c}{$k=2$} & \multicolumn{2}{|c|}{$k=4$}\Tstrut \\
     & error $E_{\tau}$ & EOC & error $E_{\tau}$ & EOC\Bstrut \\ \hline
    $2.50 \cdot 10^{-01}$ & $7.65 \cdot 10^{-05}$ & - & $2.03 \cdot 10^{-05}$ & -\Tstrut \\
    $1.28 \cdot 10^{-01}$ & $3.33 \cdot 10^{-05}$ & $1.25$ & $2.55 \cdot 10^{-06}$ & $3.11$ \\
    $6.49 \cdot 10^{-02}$ & $6.53 \cdot 10^{-06}$ & $2.39$ & $1.37 \cdot 10^{-06}$ & $0.91$ \\
    $3.27 \cdot 10^{-02}$ & $1.15 \cdot 10^{-06}$ & $2.53$ & $1.20 \cdot 10^{-07}$ & $3.55$ \\
    $1.64 \cdot 10^{-02}$ & $1.03 \cdot 10^{-08}$ & $6.83$ & $2.95 \cdot 10^{-09}$ & $5.37$ \\
    $8.21 \cdot 10^{-03}$ & $7.48 \cdot 10^{-11}$ & $7.13$ & $7.99 \cdot 10^{-12}$ & $8.55$ \\
    $4.11 \cdot 10^{-03}$ & $4.67 \cdot 10^{-12}$ & $4.00$ & $2.49 \cdot 10^{-11}$ & - \\
    $2.06 \cdot 10^{-03}$ & $1.31 \cdot 10^{-11}$ & - & $8.12 \cdot 10^{-12}$ & $1.62$ \\
    $1.03 \cdot 10^{-03}$ & $2.47 \cdot 10^{-11}$ & - & $2.18 \cdot 10^{-12}$ & $1.90$\Bstrut \\
\hline
      & & $\hspace{-11pt}\diameter\; 2.72$ & & $\hspace{-11pt}\diameter\; 2.92$\Tstrut\Bstrut \\ \hline
\end{tabular}
}
		\caption{$\zs$}
	\end{subtable}
	\hfill
	\caption{
		Nodal superconvergence in $\tau$ for several polynomial degrees $k$, for $\quadnodes = k$, and $\projnodes = 2k$, for the space-discrete (fixed mesh size $h = 1$) regularized ($\eps=10^{-10}$) porous medium equation for $q=2$. 
	}
	\label{fig:porous_medium_q2.0_eps1e-10_varying_degree_different_sampling}
\end{table}

\begin{table}
	\hfill
	\begin{subtable}[t]{.49\textwidth}
		\centering
		\scalebox{0.8}{
\begin{tabular}{|c|cc|cc|}
\hline
    \multirow{2}{\widthof{$\tau$}}{$\tau$} & \multicolumn{2}{|c}{$k=2$} & \multicolumn{2}{|c|}{$k=4$}\Tstrut \\
     & error $E_{\tau}$ & EOC & error $E_{\tau}$ & EOC\Bstrut \\ \hline
    $2.50 \cdot 10^{-01}$ & $1.07 \cdot 10^{-01}$ & - & $4.69 \cdot 10^{-02}$ & -\Tstrut \\
    $1.28 \cdot 10^{-01}$ & $8.36 \cdot 10^{-02}$ & $0.36$ & $1.79 \cdot 10^{-01}$ & - \\
    $6.49 \cdot 10^{-02}$ & $7.23 \cdot 10^{-03}$ & $3.60$ & $5.82 \cdot 10^{-02}$ & $1.66$ \\
    $3.27 \cdot 10^{-02}$ & $4.20 \cdot 10^{-02}$ & - & $2.05 \cdot 10^{-02}$ & $1.52$ \\
    $1.64 \cdot 10^{-02}$ & $2.02 \cdot 10^{-02}$ & $1.06$ & $1.32 \cdot 10^{-02}$ & $0.64$ \\
    $8.21 \cdot 10^{-03}$ & $1.11 \cdot 10^{-02}$ & $0.86$ & $1.86 \cdot 10^{-02}$ & - \\
    $4.11 \cdot 10^{-03}$ & $2.15 \cdot 10^{-02}$ & - & $5.27 \cdot 10^{-03}$ & $1.82$ \\
    $2.06 \cdot 10^{-03}$ & $5.77 \cdot 10^{-03}$ & $1.90$ & $8.04 \cdot 10^{-03}$ & - \\
    $1.03 \cdot 10^{-03}$ & $6.20 \cdot 10^{-03}$ & - & $5.44 \cdot 10^{-03}$ & $0.56$\Bstrut \\
\hline
      & & $\hspace{-11pt}\diameter\; 0.52$ & & $\hspace{-11pt}\diameter\; 0.39$\Tstrut\Bstrut \\ \hline
\end{tabular}
}
		\caption{$z_{\text{B}}$}
	\end{subtable}
	\hfill
	\begin{subtable}[t]{.49\textwidth}
		\centering
		\scalebox{0.8}{
\begin{tabular}{|c|cc|cc|}
\hline
    \multirow{2}{\widthof{$\tau$}}{$\tau$} & \multicolumn{2}{|c}{$k=2$} & \multicolumn{2}{|c|}{$k=4$}\Tstrut \\
     & error $E_{\tau}$ & EOC & error $E_{\tau}$ & EOC\Bstrut \\ \hline
    $2.50 \cdot 10^{-01}$ & $1.50 \cdot 10^{-04}$ & - & $3.02 \cdot 10^{-08}$ & -\Tstrut \\
    $1.28 \cdot 10^{-01}$ & $9.76 \cdot 10^{-06}$ & $4.10$ & $1.53 \cdot 10^{-10}$ & $7.92$ \\
    $6.49 \cdot 10^{-02}$ & $6.14 \cdot 10^{-07}$ & $4.07$ & $6.39 \cdot 10^{-13}$ & $8.05$ \\
    $3.27 \cdot 10^{-02}$ & $3.83 \cdot 10^{-08}$ & $4.04$ & $2.91 \cdot 10^{-15}$ & $7.85$ \\
    $1.64 \cdot 10^{-02}$ & $2.39 \cdot 10^{-09}$ & $4.02$ & $1.19 \cdot 10^{-15}$ & - \\
    $8.21 \cdot 10^{-03}$ & $1.49 \cdot 10^{-10}$ & $4.01$ & $2.20 \cdot 10^{-15}$ & - \\
    $4.11 \cdot 10^{-03}$ & $9.34 \cdot 10^{-12}$ & $4.00$ & $2.45 \cdot 10^{-15}$ & - \\
    $2.06 \cdot 10^{-03}$ & $5.84 \cdot 10^{-13}$ & $4.00$ & $2.74 \cdot 10^{-15}$ & - \\
    $1.03 \cdot 10^{-03}$ & $3.67 \cdot 10^{-14}$ & $3.99$ & $5.08 \cdot 10^{-15}$ & -\Bstrut \\
\hline
      & & $\hspace{-11pt}\diameter\; 4.03$ & & $\hspace{-11pt}\diameter\; 7.94$\Tstrut\Bstrut \\ \hline
\end{tabular}
}
		\caption{$\zs$}
	\end{subtable}
	\hfill
	\caption{
		Nodal superconvergence in $\tau$ for several polynomial degrees $k$, for $\quadnodes = k$, and $\projnodes = 2k$, for the space-discrete (fixed mesh size $h = 1$) regularized ($\eps=10^{-10}$) porous medium equation for $q=3$. 
		}
	\label{fig:porous_medium_q3.0_eps1e-10_varying_degree_different_sampling}
\end{table}


\captionsetup[subfigure]{labelfont=rm,skip = -0.1em}

\begin{figure}[h]
\hfill
	\begin{subfigure}[t]{.3\textwidth}
		\centering
		\scalebox{0.31}{\input{figures/porous_medium_q1.5_eps1e-10_barenblatt_varying_discretization.pgf}}
		\caption{$q=1.5$}
	\end{subfigure}
	\hfill
	\begin{subfigure}[t]{.3\textwidth}
		\centering
		\scalebox{0.31}{\input{figures/porous_medium_q2.0_eps1e-10_barenblatt_varying_discretization.pgf}}
		\caption{$q=2$}
	\end{subfigure}
	\hfill
	\begin{subfigure}[t]{.3\textwidth}
		\centering
		\scalebox{0.31}{\input{figures/porous_medium_q3.0_eps1e-10_barenblatt_varying_discretization.pgf}}
		\caption{$q=3$}
	\end{subfigure}
	\hfill
	\caption{
		Convergence in $\tau $ for several space-discretization parameters for the regularized ($\eps=10^{-10}$) porous medium equation with $q\in\{1.5, 2, 3\}$, for $k=\quadnodes = 4$, $\projnodes = 2k = 8$ and for the solution $z_{\text{B}}$. 
	}
	\label{fig:porous_medium_barenblatt_varying_discretization}
\end{figure}

\begin{figure}[h]
\hfill
	\begin{subfigure}[t]{.3\textwidth}
		\centering
		\scalebox{0.31}{\input{figures/porous_medium_q1.5_eps1e-10_smooth_varying_discretization.pgf}}
		\caption{$q=1.5$}
	\end{subfigure}
	\hfill
	\begin{subfigure}[t]{.3\textwidth}
		\centering
		\scalebox{0.31}{\input{figures/porous_medium_q2.0_eps1e-10_smooth_varying_discretization.pgf}}
		\caption{$q=2$}
	\end{subfigure}
	\hfill
	\begin{subfigure}[t]{.3\textwidth}
		\centering
		\scalebox{0.31}{\input{figures/porous_medium_q3.0_eps1e-10_smooth_varying_discretization.pgf}}
		\caption{$q=3$}
	\end{subfigure}
	\hfill
	\caption{
			Convergence in $\tau $ for several space-discretization parameters for the regularized ($\eps=10^{-10}$) porous medium equation with $q\in\{1.5, 2, 3\}$, for $k=\quadnodes = 4$, $\projnodes = 2k = 8$ and for the solution $\zs$. 
	}
	\label{fig:porous_medium_smooth_varying_discretization}
\end{figure}


\begin{figure}[h]
\hfill
	\begin{subfigure}[t]{.3\textwidth}
		\centering
		\scalebox{0.31}{\input{figures/porous_medium_q1.5_eps1e-10_barenblatt_energybalance.pgf}}
		\caption{$q=1.5$}
	\end{subfigure}
	\hfill
	\begin{subfigure}[t]{.3\textwidth}
		\centering
		\scalebox{0.31}{\input{figures/porous_medium_q2.0_eps1e-10_barenblatt_energybalance.pgf}}
		\caption{$q=2$}
	\end{subfigure}
	\hfill
	\begin{subfigure}[t]{.3\textwidth}
		\centering
		\scalebox{0.31}{\input{figures/porous_medium_q3.0_eps1e-10_barenblatt_energybalance.pgf}}
		\caption{$q=3$}
			\label{fig:porous_medium_barenblatt_energybalance-q3}
	\end{subfigure}
	\hfill
	\caption{
		Relative error in the energy balance with space-discrete Hamiltonian $\mathcal{H}_h$ ($h = 1$) for the regularized ($\eps=10^{-10}$) porous medium equation with $\tau=10^{-2}$ and the Barenblatt solution~$z_{\text{B}}$ as initial datum. 			
	}
	\label{fig:porous_medium_barenblatt_energybalance}
\end{figure}

\begin{figure}[h]
\hfill
	\begin{subfigure}[t]{.3\textwidth}
		\centering
		\scalebox{0.31}{\input{figures/porous_medium_q1.5_eps1e-10_smooth_energybalance.pgf}}
		\caption{$q=1.5$}
	\end{subfigure}
	\hfill
	\begin{subfigure}[t]{.3\textwidth}
		\centering
		\scalebox{0.31}{\input{figures/porous_medium_q2.0_eps1e-10_smooth_energybalance.pgf}}
		\caption{$q=2$}
	\end{subfigure}
	\hfill
	\begin{subfigure}[t]{.3\textwidth}
		\centering
		\scalebox{0.31}{\input{figures/porous_medium_q3.0_eps1e-10_smooth_energybalance.pgf}}
		\caption{$q=3$}
			\label{fig:porous_medium_smooth_energybalance-q3}
	\end{subfigure}
	\hfill
	\caption{
		Relative error in the energy balance with space-discrete Hamiltonian $\mathcal{H}_h$ ($h = 1$) for the regularized ($\eps=10^{-10}$) porous medium equation with $\tau=10^{-2}$ and the smooth solution $\zs$ as initial datum.			
	}
	\label{fig:porous_medium_smooth_energybalance}
\end{figure}

\printbibliography

\end{document}